\documentclass[preprint,12pt]{elsarticle}

\usepackage{graphicx}
\usepackage[english]{babel}  % hyphenation
\usepackage{amsmath,amsfonts,amssymb} % the usual ams stuff
\usepackage[lined,ruled]{algorithm2e}
\usepackage{bm}

\usepackage[active]{srcltx}  % inverse and forward search
\usepackage[notcite,notref]{showkeys}  % show labels in document
\usepackage{color}           % colors for highliting changes etc
\journal{Comput. Methods Appl. Mech. Engrg.}

\newcounter{rmrk}[section]

\newcommand{\set}[1]{\left\{#1\right\}}
\newcommand{\eigv}{\ensuremath{\boldsymbol{\upsilon}}}

\newcommand{\argmax}{\mathop{\mathrm{arg}\,\mathrm{max}}}
\newcommand{\argmin}{\mathop{\mathrm{arg}\,\mathrm{min}}}

\graphicspath{{figures/}}

\begin{document}

\begin{frontmatter}

%% Title, authors and addresses

%% use the tnoteref command within \title for footnotes;
%% use the tnotetext command for the associated footnote;
%% use the fnref command within \author or \address for footnotes;
%% use the fntext command for the associated footnote;
%% use the corref command within \author for corresponding author footnotes;
%% use the cortext command for the associated footnote;
%% use the ead command for the email address,
%% and the form \ead[url] for the home page:
%%
%% \title{Title\tnoteref{label1}}
%% \tnotetext[label1]{}
%% \author{Name\corref{cor1}\fnref{label2}}
%% \ead{email address}
%% \ead[url]{home page}
%% \fntext[label2]{}
%% \cortext[cor1]{}
%% \address{Address\fnref{label3}}
%% \fntext[label3]{}

\title{Quantifying Uncertainty in Material Damage from Vibrational Data}

%% use optional labels to link authors explicitly to addresses:
%% \author[label1,label2]{<author name>}
%% \address[label1]{<address>}
%% \address[label2]{<address>}

\author[tb]{T.~Butler}
\ead{butler.troy.d@gmail.com}

\author[ah]{A. Huhtala}
\ead{antti.huhtala@aalto.fi}

\author[ah]{M. Juntunen}
\ead{mojuntun@gmail.com}

\address[tb]{Department of Mathematical \& Statistical Sciences, University of Colorado Denver, Denver, CO, USA. T.~Butler's work is supported in part by the National Science Foundation (DMS-1228206) and Department of Energy (DE-SC0009279)}

\address[ah]{Aalto University School of Science, Department of Mathematics and Systems Analysis, Espoo, Finland}

\begin{abstract}
The response of a vibrating beam to a force depends on many physical parameters including those determined by material properties. Damage caused by fatigue or cracks result in local reductions in stiffness parameters and may drastically alter the response of the beam. Data obtained from the vibrating beam are often subject to uncertainties and/or errors typically modeled using probability densities. The goal of this paper is to estimate and quantify the uncertainty in damage modeled as a local reduction in stiffness using uncertain data. We present various frameworks and methods for solving this parameter determination problem. We also describe a mathematical analysis to determine and compute useful output data for each method. We apply the various methods in a specified sequence that allows us to interface the various inputs and outputs of these methods in order to enhance the inferences drawn from the numerical results obtained from each method. Numerical results are presented using both simulated and experimentally obtained data from physically damaged beams.
\end{abstract}

\begin{keyword}
%% keywords here, in the form: keyword \sep keyword
ensemble Kalman filter \sep Tikhonov regularization \sep inverse problems \sep parameter estimation \sep damage identification \sep measure theory 

%% MSC codes here, in the form: \MSC code \sep code
%% or \MSC[2008] code \sep code (2000 is the default)
\end{keyword}

\end{frontmatter}

%%
%% Start line numbering here if you want
%%
% \linenumbers

%% main text

% ----------------------------------------------------------------
% ----------------------------------------------------------------
\section{Introduction}
% ----------------------------------------------------------------
% ----------------------------------------------------------------

%{\bf Will fill this in later by moving content from other sections. Much of the descriptions of the methods on a high level can go here - see Section~\ref{S:Methods} discussion.}

A problem of significant importance in the field of uncertainty quantification (UQ) is parameter identification in a computational model given uncertain (i.e.~noisy) measurements. At its core, this is an inverse problem and many methods have been developed to solve it under certain assumptions. In this paper, we present three different UQ methodologies developed and analyzed with respect to distinct modeling frameworks and assumptions. We consider how the different frameworks and methodologies may be used in a complimentary manner. The goal is to not only improve the quantitative results, analysis, and inferences of any single method, but also to provide a more complete description and quantification of the uncertainty in the parameter estimate.

\begin{figure}[htb]
\centering
\includegraphics[trim=0mm 0mm 0mm 0mm, clip, scale=0.45]{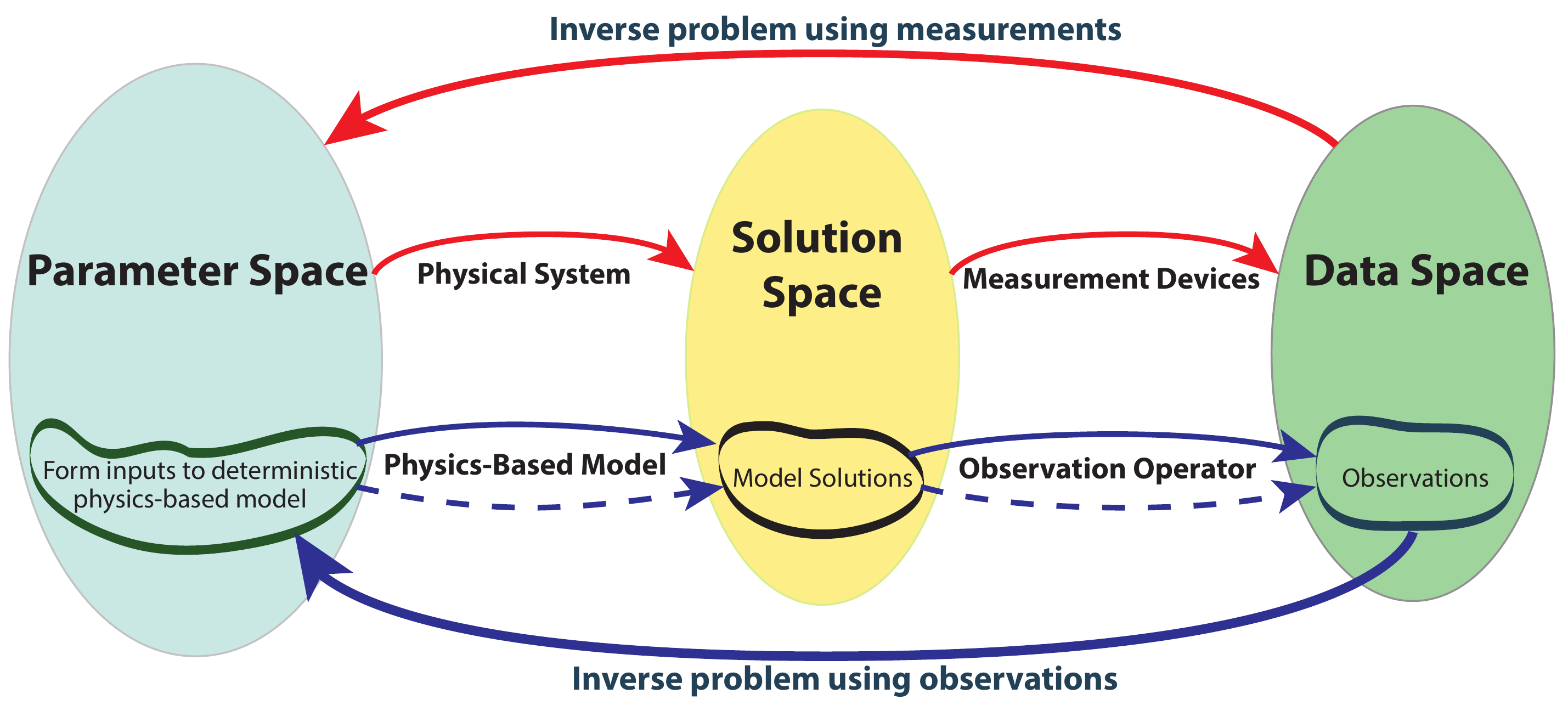}
\caption{The parameter space contains all the inputs to the physical system including the subspace forming inputs to the deterministic physics-based model. The physical system defines a map to the solution space of all possible responses of the system including the subspace of model solutions defined by the solution operator of the model. The data space contains all the possible measurement data mapped to by measurement device recordings of model solutions including those we can mathematically model as an observation operator on the subspace of model solutions. The dotted arrows denote the numerical approximation to the model and observation maps defined by a consistent numerical scheme such that the numerical approximations converge to the exact model solutions and observations assuming exact arithmetic operations.}
\label{f:Abstract_flowchart}
\end{figure}

A basic flowchart useful for defining the modeling frameworks and assumptions of the methods considered here is shown in Figure~\ref{f:Abstract_flowchart}. Ideally, we want to solve the inverse problem using the set of all physically possible measurements denoted by the red arrows directing a closed system at the top of Figure~\ref{f:Abstract_flowchart}. In this inverse problem, the solution is completely observable in space-time, and the goal is to determine all parameters, environmental effects, initial and boundary conditions, etc. Except in the most trivially controlled systems, this is an intractable problem. We instead concern ourselves with the inverse problem using observations from the model denoted by the blue arrows directing a closed system at the bottom of Figure~\ref{f:Abstract_flowchart}. The observations we choose to invert are often informed by measurements and exist in some finite dimensional space. It is not uncommon that the number of observations is orders of magnitude less than the number of unknown parameters. Prior knowledge of the system and its components is often used to define a domain within the subspace of all inputs to the model.
%For example, knowing whether or not a material being studied is rubber or polystyrene can change the domain of possible Young's modulus values from $[0.01,0.1]$ to $[3,3.5]$.
Applying some knowledge to restricting the domain of model parameters depicted in Figure~\ref{f:Abstract_flowchart} is crucial given the constraint of finite computational resources, e.g. we expect the search for a Young's modulus along the half-line $[0,\infty)$ to be computationally inefficient compared to a search within a bounded subinterval.

Another consideration in the formulation and analysis of UQ methods is that the physical model and observation map are discretized resulting in approximate maps from parameters to solutions and from solutions to observations (denoted by dotted lines in Figure~\ref{f:Abstract_flowchart}). Assuming consistent numerical approximations are implemented, these maps converge to the exact model and observation maps as the discretization parameters are refined. Thus, in the limit, the inverse problem can be formulated equivalently with respect to either the exact or computational maps between the spaces. In the numerical results, we use a fixed discretization of the space-time domain to define the numerical model.% and discuss the affects of this choice on the results.

The first UQ method we consider is the well known Ensemble Kalman Filter (EnKF) \cite{EvensenBook}. Abstractly, the model is used to ``evolve'' an ensemble of sample parameters forward in time while (often unknown) models of both process and measurement noise are used to form an invertible operator with a well-defined minimum variance solution in the parameter domain given the observations.  In Figure~\ref{f:Abstract_flowchart}, random perturbations are added in both the subspaces of model solutions and observations. In other words, the arrows representing the solution operator to the physical model and observation operator are essentially replaced by statistical maps between the  spaces. Any discrepancies between the physical solution and the model solution are accounted for in the data space. These discrepancies are then weighted by the so-called Kalman gain, and a parameter estimate is subsequently determined. This process is repeated in time as new data become available. The EnKF is typically straightforward to implement and naturally provides a quantification of uncertainty for its point estimate (i.e.~the mean) at any time via a covariance. Of all UQ methods, it is perhaps the easiest to implement in so-called ``real-time'' since at least algorithmically it is straightforward to use data as they become available. However, there is often no guarantee that the solution is physically meaningful as physical relationships between parameters, model solutions, and observations are essentially discarded in favor of statistical relationships. Also, the sequential nature of the updates via noisy data often means that the ``solution'' to the problem may never converge to a local or global solution. Nevertheless, the EnKF has proven useful in practice and it may be used to enhance our belief in the plausibility of a solution obtained via other UQ methods. We specifically employ the EnKF in order to provide a set of reasonable constraints on initial conditions and parameter domains in the optimization algorithms used to solve the problem via the second UQ method we consider: regularization.

Regularization is perhaps the most common approach for solving deterministic inverse problems. The solution operator from parameter space to solution space and/or observation space is ``regularized'' by adding an invertible operator (often called a ``penalty term'') in a linear case (e.g.~as done in Tikhonov regularization). The result is in essence a new model that has a well-posed inverse in the parameter domain (i.e.~we replace the composition of maps defined by the bottom left-to-right arrows in Figure~\ref{f:Abstract_flowchart} by a new regularized map). One may alternatively take a Bayesian point of view where the solution to the regularized model can be interpreted as the maximum a posteriori estimate to a posterior density, which assumes the model is actually statistical instead of deterministic. In either approach, when the map from parameter domain to observations is nonlinear, which is often the case even for a linear model with respect to the state vector, the regularized problem is essentially one of a nonlinear least squares problem solved by optimization, e.g.~using minimization techniques.  Similar to the EnKF, we note that the model is fundamentally changed. Any discussion of the ``closeness'' of the identified ``optimal'' parameter to an actual solution of the parameter identification problem may only make some sense within the asymptotic region of ``vanishing regularity.'' While altering the model operator may represent a serious loss of information about the model behavior, the method of regularization has nonetheless proven quite useful in practice for identifying input parameters that produce model outputs close to observed data. A question that must be addressed when using regularization is the appropriateness of the ``penalty'' term used to enforce some prior knowledge of the structure of the parameter. In other words, what prior evidence is available to suggest certain configurations or spatially variable structures of parameters are preferable to others? Another issue is that convergence of the optimization method to a global solution is often difficult to guarantee and can be strongly influenced by the initial guess. The EnKF may prove useful in providing more informed initial guesses while also providing numerical evidence that validates the type/structure of penalty imposed.

The third method we consider is based on a computational measure-theoretic approach to inverse problems for deterministic models \cite{Breidt2011,Butler2012a,Butler2013a, Butler2014a, Butler2014b}. This method fully exploits the geometric structures imposed on the parameter space by the physics defining the model operator when inverting a probability measure on the observation space to input parameters. The approach for inverting observations can be described as carrying out analysis in the space of set-valued inverses (i.e.~generalized contours) defining the solution to the physical inverse problem. This analysis on the generalized contours can be made concrete by explicit approximation for individual observations, see \cite{Breidt2011,Butler2012a}. The point of view for this method is that the model is not considered ill-posed with respect to the generalized contours, so these generalized contours are dealt with directly. In \cite{Butler2014a}, the additional geometric complications that arise from inverting multiple observations are handled using an altered approach based on approximating events rather than the manifolds defining the generalized contours. This method inverts a probability measure on observations to a uniquely defined probability measure on both the generalized contours and original input parameter domain. In other words, a probability measure is obtained on the parameter domain that propagates through the physical model to the exact probability measure imposed on the observations. Explicitly, we formulate and solve the inverse problem with respect to the Sigma-Algebras on COntour Maps (SACOM). While this method is designed explicitly for problems such as the one considered here, there are a number of practical issues to consider including defining the input parameter domain and the approximations of events in this domain. Since we compute a probability measure on the entire parameter domain, ensuring a uniform accuracy in the approximate measure and/or its density can be computationally expensive especially when the parameter domain is high-dimensional or parameters vary across several orders of magnitude. We exploit the results from regularization to define restricted parameter domains of likely solutions to the inverse problem on which the probability measure is directly computed and analyzed so that computational resources are used most efficiently.

To make these ideas less abstract, we apply the various UQ methods to a vibrating beam modeled under the Euler-Bernoulli beam theory.
% where it is assumed that there are infinitesimal rotations and strains, and the beam is in pure bending.
We observe a vibrating beam's response in order to determine the presence of damage such as cracks or fatigue that generally affect material properties (in this case by a local reduction in stiffness) while preserving mass. Observations are obtained via accelerometers placed at various known locations along the beam. The observations contain uncertainty due to measurement errors. We note that the particular model studied here is a well-defined deterministic problem since the model is defined by a physics-based deterministic equation and {\em any uncertainty} from the point of view of the model is due to either uncertain input parameters, i.e.~stiffness coefficients, and/or uncertain data, e.g.~due to measurement noise.

The outline of the paper is as follows. In Sec.~\ref{S:Model}, we derive a completely discretized model for a vibrating beam with localized mass-preserving damage, modeled as local reductions in stiffness (e.g.~as results from cracks or fatigue). We then describe our discretized observation model. In Sec.~\ref{S:Methods}, we provide an overview of the methodologies used to estimate and quantify the uncertainty for the local stiffness parameters. In Sec.~\ref{S:numerical_framework}, we describe the basic computational framework for the numerical results and specifically describe how to pre-process the data sets to provide useful output data for each of the UQ methods.  In Sec.~\ref{S:Numerics}, we provide numerical results for various damage cases using both simulated and experimental data.  Concluding remarks follow in Sec.~\ref{S:Conclusions}.

% ----------------------------------------------------------------
% ----------------------------------------------------------------
\section{A mathematical model for a vibrating beam}\label{S:Model}
% ----------------------------------------------------------------
% ----------------------------------------------------------------
We derive the variational problem of a bending beam
%in equilibrium and appeal to Newton's second and third laws to obtain the natural extension of this problem for a bending beam
in a state of acceleration (i.e.~a vibrating beam). The differential equation for the displacement of the centerline $u(t,x)$ of the beam without external loading is
\begin{equation}
\label{eq:diff_eq}
\mu \ddot{u}(t,x) + (EIu''(t,x))'' = 0, \quad x\in(0,L), \  t \in [0,T],
\end{equation}
where $\mu$ is the mass density along the length of the beam, $E$ is the Young's modulus, and $I$ is the second moment of inertia. The quantity $EI$ is called the flexural rigidity. Above, we used the common notation $\dot{u}(t,x) = \frac{\partial}{\partial t}u(t,x)$ and $u'(t,x) = \frac{\partial}{\partial x}u(t,x)$.

We now specify the associated boundary conditions to Eq.~\eqref{eq:diff_eq}. We assume that the left end at $x=0$ is clamped and that the right end at $x=L$ is free. At the clamped end, the displacement and the rotation vanish, that is,
\begin{equation}
\label{eq:bc1}
u(t,0) = 0 \quad \text{and} \quad u'(t,0) = 0,
\end{equation}
and at the free end the moment and the shear force vanish, that is,
\begin{equation}
\label{eq:bc2}
EI u''(t,L) = 0 \quad \text{and} \quad EIu'''(t,L) = 0.
\end{equation}

To complete the model definition, we need to know the initial conditions, i.e.~the displacement $u(0,x)$ and the velocity $\dot{u}(0,x)$.% We describe the model's initial conditions in the numerical examples.

%We summarize notation used in the derivation of the variational form of the vibrating beam in Table~\ref{tab:notation} and illustrate in Figure~\ref{fig:XX}.
%
%\begin{table}[h]
%
%\begin{tabular}{ | l | l || l | l |}
%\hline
%      $L$ & Beam length &   $y_c$ & Central axis \\
%          \hline
%           $q(x)$ & Load density &  $\theta(x)$ & Angle from $y_c$ \\
%          \hline
%         $v(x)$ & Vertical displacement & $Q$ & Shear force   \\
%          \hline
%           $M$ & Bending moment & $\epsilon_{xx}$ & Horizontal strain \\
%          \hline
%          $E$ & Young's modulus & $\epsilon_{yy}$ & Vertical strain \\
%          \hline
%          $\sigma_{xx}$ & Stress in beam & $\epsilon_{xy}$ & Shear strain \\
%          \hline
%          $W_{ext}(v)$ & Work against external forcing & $W_{int}(v)$ & Potential energy  \\
%          \hline
%           $\kappa_0$, $\kappa_L$ & Torsion spring stiffness & $I$ & Second moment of area \\
%          \hline
%          $k_0$, $k_L$ & Spring stiffness &  $\sigma_x$ & Stress in beam \\
%\hline
%\end{tabular}
%\caption{Summary of notation in derivation of variational form of beam equation.}
%\label{tab:notation}
%\end{table}
%
%{\bf A figure that combines figures like 2.2 and 2.3 from Antti's thesis should go here}

\subsection{Spatial discretization} \label{S:fem}
We use a finite element method in the spatial dimension. Below, we define the solution space and the weak form of the differential equation given by Eqs.~\eqref{eq:diff_eq}--\eqref{eq:bc2}.

We use $N+1$ mesh points $0=x_1<x_2<\cdots<x_{N+1}=L$ uniformly spaced along the beam to define $N$ elements $T_i = (x_i, x_{i+1})$ for $i=1,\ldots,N$, where $x_i=(i-1)h$ and $h=L/N$. We use the standard third-order $C^1$-continuous Galerkin finite elements \cite{Argyris, Johnson} to build the solution space of continuously differentiable functions on $[0,L]$ given by,
\begin{equation}
\label{eq:Vh}
V_h := \left\{ \phi\in H^2(0,L) \  : \  \phi|_{x=0}=\phi'|_{x=0}=0, \  \phi|_{T_i}\in\mathcal{P}_3(T_i) \text{ for } i=1,\ldots, N \right\}.
\end{equation}
Here $\mathcal{P}_3(T_i)$ denotes the space of  third order polynomials on $T_i$. The essential boundary conditions at the clamped end are defined in the space $V_h$.

The discretized variational form of Eqs.~\eqref{eq:diff_eq}--\eqref{eq:bc2} is: Find $u_h(t,x) \in V_h$ for all $t>0$ such that
\begin{equation}
\label{eq:FEM_beam}
m(\ddot{u}_h,w_h)+k(u_h,w_h) = 0, \ \forall w_h\in V_h,
\end{equation}
where
\begin{align}
m(\ddot{v},w) &= \int_0^L \mu \ddot{v} w \, dx,
\label{eq:m} \\
k(v,w) &= \int_0^L EIv''w''\, dx.
\label{eq:k}
\end{align}

The space $V_h$ has $2N$ basis functions denoted by $\phi_i(x)$ and the spatially discretized solution has the form $u_h(t,x) = \sum_{i=1}^{2N} \mathbf{u}_i(t) \phi_i(x)$. Here, $\mathbf{u}(t)\in\mathbb{R}^{2N}$ is a vector where the components $\mathbf{u}_i(t) \in \mathbb{R}$ are the coefficients of the basis functions at time $t$. If $u_h(t,x)$ solves Eq.~\eqref{eq:FEM_beam} for any $\phi_j$, $j=1,\ldots,2N$, then $u_h(t,x)$ solves Eq.~\eqref{eq:FEM_beam} for any $w_h \in V_h$. Thus we transform the problem of solving Eq.~\eqref{eq:FEM_beam} into the second-order linear system
\begin{equation}
\label{eq:ode_beam}
M\ddot{\mathbf{u}} + {K}\mathbf{u} = 0,
\end{equation}
where the mass and stiffness matrices are
\[
\ [M]_{ji} = m(\phi_i,\phi_j) \quad \text{and} \quad [K]_{ji}=k(\phi_i,\phi_j).
\]

\subsection{Damage model} \label{S:damage}
We are interested in observing a vibrating beam's response to the presence of damage. We assume that the damage, such as cracks or fatigue, reduces the flexural rigidity while preserving mass. The damage is modeled as a local element-wise reduction in stiffness. This model is not an exact representation of a crack in a beam, as the element size is typically much larger than the crack size, and no stress concentrations are considered near the crack. Hence the amount of damage may not be quantitatively inferred from this model, but the model is capable of localizing the damage. More elaborate damage models have been considered e.g. in \cite{friswell1,friswell2,pandey,friswell2002crack}.

The stiffness matrix $K$ can be written as a sum of stiffness matrices for each element, that is, $K=\sum_{i=1}^N K_{T_i}$ where $K_{T_i}$ is the stiffness matrix for element $T_i$. Assuming element-wise constant damage, we model the damage as a perturbation
%$\Delta K$
to the stiffness matrix $K$ such that
\begin{equation}\label{eq:stiffness_damage}
	\hat{K}(\mathbf{d}) = %K + \Delta K =
	\sum_{i=1}^N (1-\mathbf{d}_i) K_{T_i},
\end{equation}
where $\mathbf{d} \in \mathbb{R}^N$, and $\mathbf{d}_i\in [0,1)$ is the stiffness reduction in element $T_i$. Cases $\mathbf{d}_i=0$ and $\mathbf{d}_i=1$ imply no damage and complete loss of stiffness in element $T_i$, respectively. Hence, solutions for a vibrating beam that has suffered local damage approximately solve the second-order linear system of ODEs
\begin{equation}\label{eq:ode_beam_damage}
%	M\ddot{U} + (K+\Delta K)U = 0. % \mathbf{0}_{n\times 1}
	M\ddot{\mathbf{u}} + \hat{K}(\mathbf{d})\mathbf{u} = 0. % \mathbf{0}_{n\times 1}
\end{equation}

We restrict focus to the well established Euler-Bernoulli beam theory and discretize both the model and its parameters at coarse scales. The arguments for using lower or higher fidelity models bring into question the validity of the physical model and the appropriate length scales, which are beyond the scope of this paper. In other words, the physical and numerical modeling assumptions stated above create a framework where all the damage is viewed as element-wise defined loss of stiffness.

%Since we restrict focus to the well established Euler-Bernoulli beam theory and discretize both the model and its parameters at macro (i.e.~coarse) length scales, we avoid issues regarding the fidelity of the model. The arguments for using lower or higher fidelity models brings into question uncertainty in the physical model and the appropriate length scales for certain physical assumptions such as the continuum hypothesis, which is beyond the scope of this paper. In other words, the physical and numerical modeling assumptions stated above create a framework where any uncertainty in a mode shape may be viewed strictly as uncertainty in the element-wise defined stiffness parameter.

\subsection{Damping} \label{S:damping}
The vibrations in a real beam decay over time. We account for this by introducing so-called modal damping into the model. The modal damping appears as the term \cite{Argyris}
\begin{equation}
\label{eq:damping}
(\alpha M + \beta \hat{K}(\mathbf{d}))\dot{\mathbf{u}},
\end{equation}
where $\alpha \in \mathbb{R}$, $\alpha \geq 0$, and $\beta \in \mathbb{R}$, $\beta \geq 0$ are the damping parameters. The damaged ODE system Eq.~\eqref{eq:ode_beam_damage} with damping is
\begin{equation}
\label{eq:damp_beam}
M\ddot{\mathbf{u}} +(\alpha M + \beta \hat{K}(\mathbf{d}))\dot{\mathbf{u}} + \hat{K}(\mathbf{d})\mathbf{u} = 0.
\end{equation}

\subsection{Time data} \label{S:obs}

In our laboratory and numerical experiments, the measurement devices and observation operators are accelerometers attached to the beam at known locations. This implies the initial conditions are unknown since we do not have data on the initial displacements or velocities. %Because of the measurement errors, this is unfortunately true even if we knew the initial state accurately.
However, we note that the mass and stiffness matrices $M$ and $\hat{K}(d)$ are {\em independent of time} in the damage model. Thus, repeated time differentiation of Eq.~\eqref{eq:damp_beam} yields
\begin{equation*}
M{\ddot{\ddot{\mathbf{u}}}} +(\alpha M + \beta \hat{K}(\mathbf{d})){\dot{\ddot{\mathbf{u}}}} + \hat{K}(\mathbf{d})\ddot{\mathbf{u}} = 0.
\end{equation*}
Denoting the acceleration with $\ddot{\mathbf{u}}=\mathbf{a}$ gives
\begin{equation}\label{eq:compare_damp}
M{\ddot{\mathbf{a}}} +(\alpha M + \beta \hat{K}(\mathbf{d})){\dot{\mathbf{a}}} + \hat{K}(\mathbf{d})\mathbf{a} = 0.
\end{equation}
Comparing Eq.~\eqref{eq:compare_damp} to Eq.~\eqref{eq:damp_beam}, we see that the acceleration and the displacement fulfill the same equation.
%We are interested in finding the change in the flexural rigidity $d$
Since we are not interested in finding the true displacements of the system but the change in the flexural rigidity, we can interpret the acceleration data as the displacement data of the system.

%Thus we can interpret the acceleration data as displacement data

We assume the available data are defined by a relatively small number of sensors. More precisely, we assume there are $p$ sensors and $p<N$.
Let $\mathbf{y}(t)\in\mathbb{R}^p$ denote the observations at time $t$. The sensors are modeled by linear functionals of the displacement, that is, there exists
%{\color{red} Is this true with damping? Should also $\dot{U}(t)$  be included?} In our system, the measurement devices are accelerometers attached to the beam at known locations.
$B \in \mathbb{R}^{p \times 2N}$ such that
\begin{equation}\label{eq:obs_network}
\mathbf{y}(t) = B\mathbf{u}(t)
.%= BX\ddot{\eta}(t).
\end{equation}
%The last form follows from \eqref{eq:modalbase} and we note that $BX$ are the observable mode shapes.
% From eq.~\ref{eq:delta_X}, we conclude that such measurements are sensitive to localized changes in stiffness parameters. Thus, such measurements should prove useful in inverting observations for the purpose of identifying localized damage.

\subsection{Frequency data} \label{S:modal}
%Suppose we solve the $2N$ eigenvalues $\omega_i^2$ and eigenvectors $\eigv_i \in \mathbb{R}^{2N}$ of the system
For $i=1,\ldots,2N$, let $\omega_i^2$ and $\eigv_i\in\mathbb{R}^{2N}$ denote the eigenvalues and eigenvectors of the system
%Suppose we solve the $2N$ eigenvalues $\omega_i^2$ and eigenvectors $\eigv_i \in \mathbb{R}^{2N}$ of the system
\begin{equation}
\label{eq:eigen_problem}
\hat{K}(\mathbf{d}) \eigv_i = \omega_i^2 M \eigv_i.
\end{equation}
The $\omega_i$ are called the modal frequencies and the $\eigv_i$ are mode shapes of the system.
%Using a change of coordinates, eq.~(\ref{eq:ode_beam}) may be written as a system of $2N$ uncoupled harmonic oscillators with solutions
Setting $\eta_i(t)=a_i\cos(\omega_i t)+b_i\sin(\omega_i t)$, where the coefficients $a_i$ and $b_i$ depend on the initial conditions of $\mathbf{u}$ and $\dot{\mathbf{u}}$, the solutions can be written as
\begin{equation}
\label{eq:modalbase}
\mathbf{u}(t) = \sum_{i=1}^{2N} \eta_i(t)\eigv_i. % = V\mathbb{\eta}(t;\omega_i).
\end{equation}
Recalling the observation operator $B$ in Eq.~\eqref{eq:obs_network}, the observations are
\begin{equation}\label{eq:modal_obs_network}
\mathbf{y}(t) = B\mathbf{u}(t) =  \sum_{i=1}^{2N} \eta_i(t)B \eigv_i,
\end{equation}
where $B\eigv_i$ are the observable mode shapes. The mode shapes are normalized so that the Euclidean norm of the observable part is one, i.e.~$\eigv_i^TB^TB\eigv_i=1$, and that the displacement measured closest to the free end is non-negative.%, $B\eigv_i|_{x=L}\geq0$.

From Eq.~\eqref{eq:eigen_problem} it is clear that the mode shapes and the modal frequencies depend on the reduction in stiffness given by $\mathbf{d}$. When the mode shapes and the modal frequencies are sensitive to perturbations in $\mathbf{d}$, then data determined from these quantities is useful for determining the damage quantified by $\mathbf{d}$. %We show in the numerical results several scenarios of data with varying sensitivity with respect to perturbations in $\mathbf{d}$.

\subsection{Mapping parameters to data}

%In this article we consider three different frameworks and subsequent methodologies for the quantification of uncertainty of damage in the beam, namely the ensemble Kalman filter (EnKF), the Bayesian inverse, and a stochastic inverse formulated in a measure-theoretic framework. Each so-called ``inverse problem'' is formulated with respect to a specific modeling framework with different physical assumptions and using possibly different forms of the data (e.g.~using observations in time or of modal parameters).
The specification of a map from model parameters (i.e.~parameter $\mathbf{d}$) to the data as illustrated by the bottom left-to-right arrows of Figure~\ref{f:Abstract_flowchart} is not the same for each method we consider. In other words, the ``forward map'' from parameters to data requires a different specification of the ``forward solve'' of the system, depending on the method we consider.

\subsubsection{Time domain forward solve} \label{S:forward1}
For the EnKF, the forward model is based on the system \eqref{eq:damp_beam}. The forward solve requires the flexural rigidity $EI$ of the undamaged beam (assumed here to be constant), the element-wise stiffness reduction $\mathbf{d}$, and the damping parameters $\alpha$ and $\beta$. Suppose the time-interval $[0,T]$ is divided into $N_t$ time instances $0=t_0 < t_1 < \cdots < t_k < \cdots < t_{N_t} = T$. We assume the time step is constant $\delta t = t_{k+1}-t_k$. Let $\mathbf{u}_k = \mathbf{u}(t_k)$ and $\dot{\mathbf{u}}_k=\dot{\mathbf{u}}(t_k)$ denote the displacement and velocity at time $t_k$. For the numerical time discretization, the system \eqref{eq:damp_beam} is reduced to the first order system
\begin{equation*}
\begin{bmatrix} M & \\ & I \end{bmatrix}
\frac{d}{dt}{\begin{bmatrix} \dot{\mathbf{u}} \\ \mathbf{u} \end{bmatrix}}
= \begin{bmatrix} -\alpha M -\beta \hat{K}(\mathbf{d}) & -\hat{K}(\mathbf{d}) \\ I &  \end{bmatrix}
{\begin{bmatrix} \dot{\mathbf{u}} \\ \mathbf{u} \end{bmatrix}}
\end{equation*}
and discretized using the implicit mid-point rule \cite{LeVequeBook}
\begin{multline*}
\begin{bmatrix} M  + \frac{\delta t}{2} \alpha M + \frac{\delta t}{2}\beta \hat{K}(\mathbf{d}) & \frac{\delta t}{2} \hat{K}(\mathbf{d}) \\  -\frac{\delta t}{2} I& I \end{bmatrix}
\begin{bmatrix} \dot{\mathbf{u}}_{k+1} \\ \mathbf{u}_{k+1} \end{bmatrix}
\\= \begin{bmatrix} M  - \frac{\delta t}{2} \alpha M - \frac{\delta t}{2}\beta \hat{K}(\mathbf{d}) & -\frac{\delta t}{2} \hat{K}(\mathbf{d}) \\  \frac{\delta t}{2} I& I \end{bmatrix}
\begin{bmatrix} \dot{\mathbf{u}}_k \\ \mathbf{u}_k \end{bmatrix}.
\end{multline*}
The solution operator of the system above is denoted by
\begin{equation}
\label{eq:enkfSolve}
\begin{bmatrix} \dot{\mathbf{u}}_{k+1} \\ \mathbf{u}_{k+1} \end{bmatrix} = g\left( \dot{\mathbf{u}}_k, \mathbf{u}_k ;\mathbf{d},EI,\alpha,\beta \right),
\end{equation}
where we made explicit the dependence of the solution on the flexural rigidity, damage, and damping parameters. The solution also depends on, for example, the mass density $\mu$, the length of the beam $L$, the number of elements $N$, etc., but these are assumed fixed. %The starting point $(\dot{\mathbf{u}}_0,\mathbf{u}_0)$ is acquired from the initial conditions.
The observations $\mathbf{y}(t_k) = \mathbf{y}_k$ at time $t_k$ are computed with Eq.~\eqref{eq:obs_network}, i.e.
\begin{equation} \label{eq:enkfObs}
\mathbf{y}_k = B\mathbf{u}_k .
\end{equation}

\subsubsection{Frequency domain forward solve} \label{S:forward2}

For the other two methods, we use the modal basis for the forward solve. We let
\begin{equation}
\label{eq:modal_solve}
f( \mathbf{d} ) = \begin{bmatrix} \omega_1 \\ B\eigv_1 \\ \vdots \\ \omega_n \\ B \eigv_n \end{bmatrix} \in \mathbb{R}^{np+n}.
\end{equation}
denote the solution of the eigensystem described by Eq.~\eqref{eq:eigen_problem} for the $n$ smallest frequencies and the associated normalized observable modes. Here we only explicitly state the dependency on the element-wise damage $\mathbf{d}$, but as before the solution also depends on other variables, which are assumed fixed. %We explicitly stated that the solution depends on the flexural rigidity and damage, and the rest of the parameters, such as, the mass density $\mu$, the length of the beam $L$, the number of elements $N$, etc. are assumed fixed.

% ----------------------------------------------------------------
% ----------------------------------------------------------------
\section{Methods of parameter identification and uncertainty quantification}\label{S:Methods}
% ----------------------------------------------------------------
% ----------------------------------------------------------------

We provide brief descriptions and computational algorithms for the methods considered in this work. It is outside the scope of this work to consider every comparable method, variants of these methods, or all of the relevant theoretical results of the methods considered here. The EnKF and regularization methods are well established in the existing literature, and it is impossible to provide a complete list of citations to all the relevant and interesting work available. The interested reader should refer to both the pedagogical and survey citations we point out and the references contained therein.

\subsection{Data assimilation and the ensemble Kalman filter} \label{S:enkf}
% ----------------------------------------------------------------
% ----------------------------------------------------------------

%We limit the description of the widely discussed EnKF as applied to our specific application.
For a thorough pedagogical development of the EnKF, we direct the interested reader to \cite{EvensenBook, Evensen03}. A number of interesting variants of the EnKF have been developed in recent years, e.g., see \cite{Johns_and_Mandel,Pham98,WvdM00,Tippett03,ButlerJuntunen} and the references therein for a relatively small snapshot of such variants.

In the EnKF approach, we try to fit two models to the data simultaneously. The first model assumes a constant flexural rigidity without damage. The second model uses the flexural rigidity of the first model, and fits the element-wise damage $\mathbf{d}$ to the data. A seemingly easier approach would be to measure the reference flexural rigidity from the undamaged beam measurements, and then compare the results to a beam with unknown damage using a separate EnKF run. However, we found the former more effective and describe this below. This is due to environmental effects not modeled that cause changes to the dynamics of the beam between the undamaged measurements and damaged measurements.

Let the superscript ``$d0$'' denote the model which assumes zero damage. Let $\bm{m}(t_k) \in \mathbb{R}^p$ denote the displacement obtained from measurement matching the observation operator described in Section \ref{S:obs}. For our EnKF approach we duplicate it as follows
\begin{equation} \label{eq:enkfMeas}
\mathbf{\widetilde{m}}_k = \begin{bmatrix} \bm{m}(t_k) \\ \bm{m}(t_k) \end{bmatrix} \in \mathbb{R}^{2p}.
\end{equation}
%where $x_i \in \mathbb{R}$, $i=1,\ldots,p$, are the known locations of the measurement devices. Let $\mathbf{m}_k = \mathbf{m}(t_k)$ denote the physical measurements at time $t_k$.

Recall the forward solver described in Section \ref{S:forward1}. The state vector of the EnKF is
\begin{equation}
\boldsymbol{\psi}_k = \begin{bmatrix} \dot{\mathbf{u}}^{d0}_k \\ \mathbf{u}^{d0}_k  \\ \dot{\mathbf{u}}_k \\ \mathbf{u}_k  \\ EI \\ \mathbf{d} \\ \alpha \\ \beta \end{bmatrix} \in \mathbb{R}^{9N+3},
\end{equation}
where $EI \in \mathbb{R}$ is a constant flexural rigidity of the beam. The forward solve for this method is
\begin{equation} \label{eq:enkfg}
\boldsymbol{\psi}_{k+1} = \hat{g}(\boldsymbol{\psi}_{k}) =
\begin{bmatrix} g(\dot{\mathbf{u}}^{d0}_k,\mathbf{u}^{d0}_k; 0 ,EI,\alpha,\beta)  \\ g(\dot{\mathbf{u}}_k,\mathbf{u}_k; \mathbf{d} ,EI,\alpha,\beta) \\ EI \\ \mathbf{d} \\ \alpha \\ \beta \end{bmatrix},
\end{equation}
that is, we use the forward solver of Eq.~\eqref{eq:enkfSolve} for both beams and the parameters remain the same. The observations are collected using Eq.~\eqref{eq:enkfObs} and the observation operator defined to match the measurements in Eq.~\eqref{eq:enkfMeas} as
\begin{equation} \label{eq:enkfH}
\mathbf{y}_k = H \boldsymbol{\psi}_k = \begin{bmatrix} 0 & B & 0 & 0 & 0 & 0 & 0 & 0 \\ 0 & 0 & 0 & B & 0 & 0 & 0 & 0 \end{bmatrix} \begin{bmatrix} \dot{\mathbf{u}}^{d0}_k \\ \mathbf{u}^{d0}_k  \\ \dot{\mathbf{u}}_k \\ \mathbf{u}_k  \\ EI \\ \mathbf{d} \\ \alpha \\ \beta \end{bmatrix} = \begin{bmatrix} B \mathbf{u}^{d0}_k \\ B \mathbf{u}_k \end{bmatrix} \in \mathbb{R}^{2p}.
\end{equation}

Let $N_e$ denote the size of the ensemble. We denote with
\begin{equation}
X_k = \begin{bmatrix} \boldsymbol{\psi}_k^{(1)}, \ldots, \boldsymbol{\psi}_k^{(N_e)} \end{bmatrix} \in \mathbb{R}^{(9N+3) \times N_e}
\end{equation}
the collection of the ensemble members and assume the operations defined in Eq.~\eqref{eq:enkfg} and Eq.~\eqref{eq:enkfH} may also be applied column-wise to $X_k$.

The model noise is $\mathcal{E}_{\text{mod}} \in \mathbb{R}^{(9N+3) \times N_e}$, in which each column is from a normal distribution $\mathcal{N}(0,\Sigma_{\text{mod}})$ with a diagonal $\Sigma_{\text{mod}}$. Similarly, the measurement noise is $\mathcal{E}_{\text{meas}} \in \mathbb{R}^{2p \times N_e}$, in which each column is from a normal distribution $\mathcal{N}(0,\Sigma_{\text{meas}})$ with a diagonal $\Sigma_{\text{meas}}$.

Let $\mathbb{I}_{k \times l}$ denote a matrix of size $k \times l$ where each component is one and let $\text{E}(\cdot)$ denote the row-wise mean value. The EnKF method is described by Algorithm~\ref{Alg:EnKF}.

%\begin{algorithm}\caption{EnKF on $[t_k,t_{k+1}]$}\label{Alg:EnKF}
%\begin{algorithmic}
%\State Given analyzed state $X_k^a$ at $t_k$ and measurement $\mathbf{m}_{k+1}$ at $t_{k+1}$.
%\begin{enumerate}
%\item Compute forecast state: $X_{k+1} = \hat{g}(X_k^a)$.
%\item Compute perturbations to ensemble members and observations: $\hat{X}_{k+1} = X_{k+1} + \mathcal{E}_{\text{mod}}$, and $\mathbf{M}_{k+1} = \mathbf{m}_{k+1} \mathbb{I}_{1 \times N_e} + \mathcal{E}_{\text{meas}}$.
%\item Compute $A_{k+1} = \hat{X}_{k+1} - \text{E}(\hat{X}_{k+1}) \mathbb{I}_{1 \times N_e}$
%\item Compute Kalman gain: $P_{k+1} = \frac{1}{N_e-1} (H A_{k+1})^T(HA_{k+1}) + \Sigma_{\text{meas}}$.
%\item Compute analyzed state at $t_{k+1}$: $$X_{k+1}^a = \hat{X}_{k+1} + \frac{1}{N_e-1} A_{k+1} (HA_{k+1})^T P_{k+1}^{-1} ( \mathbf{M}_{k+1} - H \hat{X}_{k+1}).$$
%\end{enumerate}
%\end{algorithmic}
%\end{algorithm}

\begin{algorithm}\caption{EnKF on $[t_k,t_{k+1}]$}\label{Alg:EnKF}
 Given analyzed state $X_k^a$ at $t_k$ and measurement $\mathbf{m}_{k+1}$ at $t_{k+1}$. \\
\begin{enumerate}
\item Compute forecast state: $X_{k+1} = \hat{g}(X_k^a)$.
\item Compute perturbations to ensemble members and observations: $\hat{X}_{k+1} = X_{k+1} + \mathcal{E}_{\text{mod}}$, and $\mathbf{M}_{k+1} = \mathbf{\widetilde{m}}_{k+1} \mathbb{I}_{1 \times N_e} + \mathcal{E}_{\text{meas}}$.
\item Compute mean-free state: $A_{k+1} = \hat{X}_{k+1} - \text{E}(\hat{X}_{k+1}) \mathbb{I}_{1 \times N_e}$
\item Compute innovation covariance: $P_{k+1} = \frac{1}{N_e-1} (H A_{k+1})^T(HA_{k+1}) + \Sigma_{\text{meas}}$.
\item Compute analyzed state at $t_{k+1}$: $$X_{k+1}^a = \hat{X}_{k+1} + \frac{1}{N_e-1} A_{k+1} (HA_{k+1})^T P_{k+1}^{-1} ( \mathbf{M}_{k+1} - H \hat{X}_{k+1}).$$
\end{enumerate}
\end{algorithm}

% ----------------------------------------------------------------
% ----------------------------------------------------------------
\subsection{Regularization approach} \label{S:statinv}
% ----------------------------------------------------------------
% ----------------------------------------------------------------

For the regularization approach the measurement $\bm{m}_{\omega} \in \mathbb{R}^{np+n}$ is a set of modal parameters matching the observation operator $f(\mathbf{d})$ described in Section~\ref{S:forward2}. Section~\ref{S:modal_identification} describes how we obtain modal data from time-signal acceleration measurements. %Then, collect the observable mode shapes and the associated frequencies as a data vector $f(\mathbf{d})$ as in Eq.~\eqref{eq:modal_solve}.
%\begin{equation}
%\mathbf{m} =
%\begin{bmatrix}
%\omega_1 \\
%B \eigv_1 \\
%\vdots \\
%\omega_p \\
%B \eigv_p
%\end{bmatrix},
%\end{equation}
%where $p$ is the number of sensors.

The map from parameters to the data is given by
\begin{equation}
\mathbf{m}_{\omega} = f(\mathbf{d}) + \bm{\epsilon},
\end{equation}
where $\bm{\epsilon} \in \mathbb{R}^{np+n}$ is additive noise and $f(\bm{d})$ is the observation operator defined in Eq.~\eqref{eq:modal_solve}. The damage identification problem is to find $\mathbf{d}$ given $\mathbf{m}_\omega$.

We derive the regularization method by the use of Bayesian statistics. Thus we consider the observation, the damage parameter vector and the noise as random variables. From the model, the connection between the random variables is
\begin{equation}
\mathcal{M} = f(\mathcal{D}) + \mathcal{E}. \label{eq:randmeasmodel}
\end{equation}
We assume that the additive noise $\mathcal{E}$ follows a multivariate normal distribution $\mathcal{N}(\bm{\mu},\bm{\Sigma})$ with density function
\begin{equation}
p(\mathcal{E}=\bm{\epsilon}) = c \exp \left\{ -\frac{1}{2} \| \bm{S} (\bm{\epsilon} - \bm{\mu}) \|^2 \right\},
\end{equation}
where $\bm{S}^T \bm{S} = \bm{\Sigma}^{-1}$, $c$ is a normalization constant, and $\| \mathbf{x} \|^2 = \mathbf{x}^T\mathbf{x}$ denotes the vector $l^2$-norm. In addition, we assume that $\mathcal{E}$ is independent of $\mathcal{D}$, which allows us to estimate the {\em statistical parameters} $\bm{\mu}$ and $\bm{\Sigma}$ from data. The estimation process is described in Section \ref{S:measstatistics}.% Note that $\bm{\mu}$ and $\bm{\Sigma}$ are not in the physical parameter space shown on the left of Figure~\ref{f:Abstract_flowchart} as they enter into this formulation via the statistical framework.

Consider the conditional case of $\mathcal{D}=\mathbf{d}$. Then
\begin{equation}
\mathcal{M} = f(\mathbf{d}) + \mathcal{E},
\end{equation}
which states that conditionally $\mathcal{M}$ is just a shifted version of $\mathcal{E}$, with the probability distribution
\begin{equation}
p(\mathcal{M}=\mathbf{m}_\omega|\mathcal{D}=\mathbf{d}) = c \exp \left\{ -\frac{1}{2} \| \bm{S} (\mathbf{m}_\omega - f(\mathbf{d}) - \bm{\mu}) \|^2 \right\}.
\end{equation}
This distribution is called the likelihood of $\mathcal{D}$. % This knowledge of the conditional distribution also allows us to estimate the statistical parameters $\bm{\mu}$ and $\bm{\Sigma}$ from the measurement data. This process is described in Section \ref{S:measstatistics}.
Bayes' theorem then states that
\begin{equation}
p(\mathcal{D}=\mathbf{d}|\mathcal{M}=\mathbf{m}_\omega) = \frac{p(\mathcal{M}=\mathbf{m}_\omega|\mathcal{D}=\mathbf{d})p(\mathcal{D}=\mathbf{d})}{p(\mathcal{M}=\mathbf{m}_\omega)}.
\end{equation}
Thus, the posterior distribution, that is, the above conditional probability distribution of $\mathcal{D}$, can be expressed through the likelihood of $\mathcal{D}$ and the marginal distributions $p(\mathcal{D}=\mathbf{d})$ and $p(\mathcal{M}=\mathbf{m}_\omega)$. The distribution $p(\mathcal{D}=\mathbf{d})$ is called the prior distribution, as it expresses information about the damage vector prior to obtaining the outcome $\mathbf{m}_\omega$. The value of $p(\mathcal{M}=\mathbf{m}_\omega)$ acts as a normalizing constant.

We choose the prior distribution as a truncated multivariate normal distribution, such that
\begin{equation}
p(\mathcal{D}=\mathbf{d}) =
\begin{cases}
0 & \exists i \text{ such that } \mathbf{d}_i < 0 \\
0 & \exists i \text{ such that } \mathbf{d}_i \ge 1 \\
c \exp \{-\frac{1}{2 \lambda^2} \| \mathbf{d} \|^2 \} & \text{otherwise},
\end{cases}
\end{equation}
where $c$ is a normalization constant. The truncation gives zero probability to the disallowed damage parameters, and otherwise prefers small values of the damage parameters over large values. The parameter $\lambda$ sets the variance of the distribution. This choice of prior distribution is justified physically by the knowledge that the damage in the beam for the simulated and experimental cases is always localized. In the absence of such physical knowledge, we may validate such a choice with the analysis of the EnKF results as we show in the numerical results.

The estimate for $\mathbf{d}$ we consider is the maximum a posteriori estimate, which is the maximum of the posterior distribution
\begin{equation}
\begin{aligned}
\mathbf{d}_{MAP} & = \argmax_{\mathbf{d}} \, p(\mathcal{D}=\mathbf{d}|\mathcal{M}=\mathbf{m}_\omega) \\
& = \argmin_{0 \le \mathbf{d}_i < 1} \left\{ -\log p(\mathcal{D}=\mathbf{d}|\mathcal{M}=\mathbf{m}_\omega)\right\} \\
& = \argmin_{0 \le \mathbf{d}_i < 1} \left\{ \| \bm{S} (\mathbf{m}_\omega - f(\mathbf{d}) - \bm{\mu}) \|^2 + \frac{1}{\lambda^2} \| \mathbf{d} \|^2 \right\}.
\end{aligned}
\end{equation}
This problem is a non-linear least squares optimization problem, which is solved using the Levenberg-Marquardt algorithm \cite{marquardt, bazaraa}. For a more thorough exposition on the topics of regularization and connections to Bayesian inference, we direct the interested reader to \cite{Tarantola,Somersalo2005}.

% ----------------------------------------------------------------
% ----------------------------------------------------------------
\subsection{Measure theory and sigma-algebras on contour maps}
% ----------------------------------------------------------------
% ----------------------------------------------------------------

Below, we give an intuitive motivation for the framework, problem formulation, and solution method referred to as SACOM (sigma-algebras on contour maps). We direct the interested reader to \cite{Breidt2011, Butler2012a, Butler2013a, Butler2014a, Butler2014b} for more details about the measure-theoretic framework, theoretical results including convergence of probabilities, computational algorithms, and a full error analysis.

Given perfect noise-free data of mode shapes and frequencies, the damage identification problem is a deterministic inverse problem. On any domain describing the possible damage parameters, for a fixed computational model, we may define the solution to the inverse problem as the damage parameter for which the corresponding observations from the model match, or are closest in some norm, to the measurements.
However, it is often the case that such a solution is not unique due either to the dimension of observation space being less than the dimension of the parameter space and/or nonlinearities in the map between these spaces leading to an inherently set-valued inverse map. Typical assumptions on the maps between parameter, solution, and output data spaces\footnote{For example, the observation operator for the model is differentiable with respect to the parameters.} guarantee that the set-valued inverses are described by piecewise-smooth manifolds of dimension $N-p$ in the parameter space. We call these set-valued inverses {\em generalized contours} and the description of generalized contour events  defines a type of contour map in the parameter space.

The uncertainties in the measurements, e.g.~due to measurement noise, are often described in stochastic terms commonly in the form of a probability density. The result is a type of stochastic inverse problem for deterministic models: Given a probability density on observations, determine the corresponding probability density on the parameter space. Introducing stochasticity in this framework does not change the fact that the inverse map is set-valued, and it can be shown that there exists a unique probability density on the contour map \cite{Breidt2011, Butler2014a}. However, we are generally interested in a probability measure or density on the original parameter space not on some possibly complicated contour map. Applying the Disintegration Theorem \cite{Dellacherie_Meyer} and an ansatz exploiting the geometric information inherent in the physical map defines a unique probability measure on the original parameter domain \cite{Butler2014a}. Several computational measure-theoretic algorithms have been developed and analyzed to approximate this unique probability measure on the original parameter domain \cite{Breidt2011, Butler2014a, Butler2014b}.

We use the results of regularization to define a domain $D$ in the parameter subspace of model inputs in which to solve the inverse problem. Here, we let $M$ denote the space of observations informed by uncertain measurements on which a probability density $\rho_M$ is given. In this work, we compute probability approximations on $D$ using the non-intrusive sample based Algorithm~\ref{Alg:measure-theory} below, which uses the same forward solve $f(\bm{d})$ as the regularization method. This sample based algorithm is particularly useful for high dimensional parameter spaces \cite{Butler2014b}. For the sake of simplicity, we consider the case where the random parameter samples are independent identically distributed (i.i.d.) and sampled uniformly with respect to the underlying volume measure on $D$.

%\begin{algorithm}\caption{A Monte Carlo Approximation of SACOM}\label{Alg:measure-theory}
%\begin{algorithmic}
%\State Let $\set{\mathbf{d}^{(j)}}_{j=1}^J\subset{D}$ denote uniform i.i.d.~random samples, and $\set{\mathcal{V}_j}_{j=1}^{J}\subset D$ denote the associated Voronoi tessellation of ${D}$.
%\For{$j=1,\ldots,N$}
%\State Assign value $f^{(j)}=f(\mathbf{d}^{(j)},EI)$ to $\mathcal{V}_j$.
%\EndFor
%\State Generate partition $\set{\mathcal{I}_i}_{i=1}^I\subset M$.
%\State Fix and normalize $\rho_{\mathcal{D},I} =\sum_{i=1}^I p_i \mbf{1}_{\mathcal{I}_i}(q)$.
%\State Initialize $I\times 1$ counting vector $\mathbf{c}$ and $J\times 1$ pointer vector $\mathbf{i_o}$ to zeros.
%\For{$j=1,\ldots,J$}
%	\State Set $i=1$ and \textbf{flag}$=0$
%	\While{$i\leq I$ and \textbf{flag}$=0$}
%		\If{ $\mathbf{f}^{(j)}\in \mathcal{I}_i$} $\mathbf{c}(i) = \mathbf{c}(i)+1$, $\mathbf{i_o}(j)=i$, \textbf{flag}$=1$.
%		\Else{ $i=i+1$}.
%		\EndIf
%	\EndWhile
%	
%\EndFor
%\For{$j=1,\ldots,J$}
%	\State Set $P_{D,J}(\mathcal{V}_j)$ to $p_{\mathbf{i_o}(j)}/\mathbf{c}(\mathbf{i_o}(j))$
%\EndFor
%\end{algorithmic}
%\end{algorithm}

\definecolor{dark-gray}{gray}{0.3}

\begin{algorithm}\caption{Sample Based Approximation of the Probability Measure using SACOM}\label{Alg:measure-theory}
Let $\set{\mathbf{d}^{(j)}}_{j=1}^J$ denote a set of samples and $\set{\mathcal{V}_j}_{j=1}^{J}$ denote the associated Voronoi tessellation of $D$. \\
%\For{$j=1,\ldots,N$}
%{
%  Assign value $f^{(j)}=f(\mathbf{d}^{(j)})$ to $\mathcal{V}_j$.
%}
Generate a tessellation $\set{\mathcal{I}_i}_{i=1}^I$ of $M$. \\
%Compute $\rho_{M,I} =\sum_{i=1}^I p_i \mbf{1}_{I_i}(f)$ where $p_i=\int_{I_i} \rho_{M}(f)\, d\mu_{M}$. \\
Compute the probabilities $p_i =\int_{I_i} \rho_{M}\, \mathrm{d}\mu_{M}$ for each $I_i$.
Initialize counting vector $\mathbf{c} \in \mathbb{Z}^I$ and index vector $\mathbf{k} \in \mathbb{Z}^J$ to zeros. \\
\For{$j=1,\ldots,J$}
{
% Set $i=1$ and \textbf{flag}$=0$
% \While{$i\leq I$ and \textbf{flag}$=0$}
% {
%	\eIf{ $f(\bm{d}^{(j)})\in \mathcal{I}_i$}
%	{
%	  $\mathbf{count}(i) = \mathbf{count}(i)+1$, $\mathbf{i_o}(j)=i$, \textbf{flag}$=1$.
%    }
%	{$i=i+1$}.
% }
Find the unique $k$ such that $f(\bm{d}^{(j)}) \in \mathcal{I}_k$. \\
Accumulate $\bm{c}(k) = \bm{c}(k) + 1$. \\
Set $\bm{k}(j) = k$.
}
\For{$j=1,\ldots,J$}
{
% Set $P_{D,J}(\mathcal{V}_j)$ to $p_{\mathbf{i_o}(j)}/\mathbf{c}(\mathbf{i_o}(j))$
Compute the probability of Voronoi cell $\mathcal{V}_j$: \\
$P(\mathcal{V}_j) = p_{\mathbf{k}(j)}/\mathbf{c}(\mathbf{k}(j))$
}
\end{algorithm}

Following Algorithm~\ref{Alg:measure-theory}, we may compute a {\em probability counting measure} \cite{Butler2014b} on the parameter space such that the probability of an arbitrary event of parameters, denoted by $A$, may be computed by
\[
	P(A) = \sum_{j=1}^J P(\mathcal{V}_j)\, \chi_A(\mathbf{d}^{(j)}),
\]
where $\chi_A$ is the indicator function of set $A \subset D$.

The measurement data used in this method is the same as used in the regularization approach, which matches the observation operator $f(\bm{d})$ described in Section~\ref{S:forward2}. The probability distribution $\rho_M$ of the algorithm is defined using the same measurement noise model that was used with the regularization method, i.e.
\begin{equation}
\mathcal{M} = f(\mathcal{D}) + \mathcal{E},
\end{equation}
where $\mathcal{E}$ is $\mathcal{N}(\bm{\mu},\bm{\Sigma})$ distributed. For a given measurement $\bm{m}_\omega$, the observation then has the distribution
\begin{equation}
\rho_M(\bm{x}) = c \exp\left\{ -\frac{1}{2} \| \bm{S} ( \bm{x} - \bm{m}_\omega + \bm{\mu} ) \|^2 \right\},
\end{equation}
where $\bm{S}^T \bm{S} = \bm{\Sigma}^{-1}$ and $c$ is a normalizing constant. The statistical parameters for the distribution are estimated from the data. This procedure is described in Section~\ref{S:measstatistics}.

% ----------------------------------------------------------------
% ----------------------------------------------------------------
\section{Measurement framework}\label{S:numerical_framework}
% ----------------------------------------------------------------
% ----------------------------------------------------------------

We describe the experimental measurement setup as well as our simulated measurements. We also provide details on the computational approaches we use to process the data for use in the various UQ methodologies.

%We also summarize the basic numerical setup of the forward models defining the maps from parameter space to solution space.

\subsection{Experimental measurement data}

Experimental data were obtained from a steel cantilever beam with various levels of manufactured damage. The beam was 1400 mm long in the $x$-direction, 60 mm wide in the $z$-direction and 5 mm thick in the $y$-direction. The vibration of the beam was measured with seven accelerometers placed along its length with a sampling rate of 512 Hz. The vibrations were only measured in the $y$-direction, which was the principle direction of the vibration. To excite the vibrations of the beam, it was bent from the free end by hand, and let loose to vibrate freely.

Before any damage was induced, five initial sets of measurements were recorded. Each measurement set consists of 30 seconds of vibration after excitation. Using a hacksaw an approximately 1 mm wide (in the $x$-direction) and 5 mm deep (in the $z$-direction) slot was cut through the beam at 260 mm from the fixed end. Five sets of measurements of 30 seconds were again recorded from the now damaged beam after excitation. To increase the damage further, the slot depth was subsequently increased to 10 mm, 15 mm, and finally 20 mm. Five 30 second measurement sets were recorded for each of these slot depths. We refer to these measurement sets in terms of their associated damage cases 1, 2, 3 and 4, respectively.

\begin{figure}[htb]
\begin{center}
\includegraphics[scale=0.65]{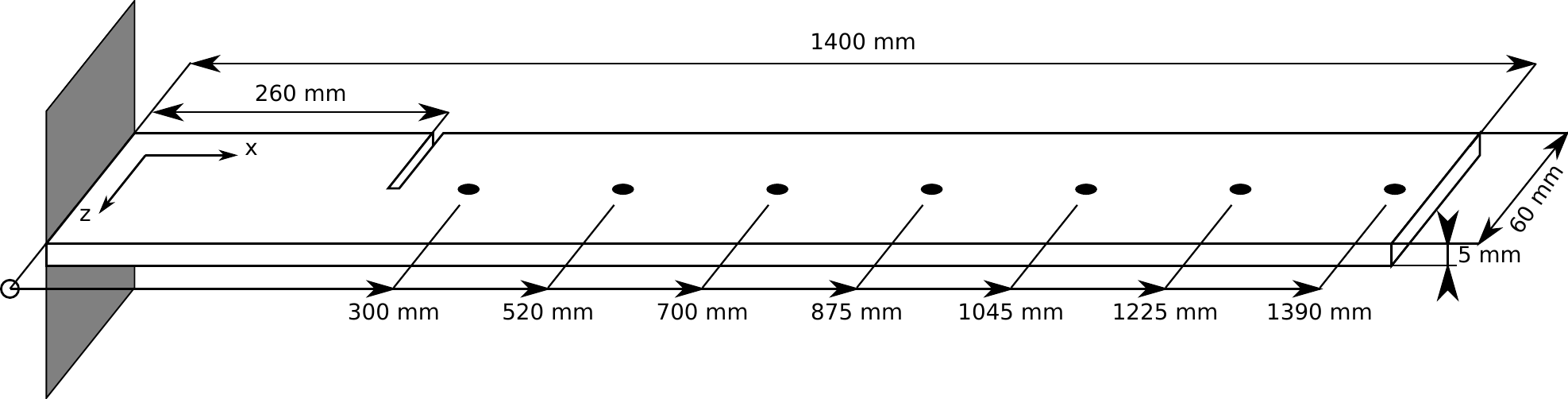}
\caption{Experiment setup}
\label{F:experiment}
\end{center}
\end{figure}

\begin{figure}[htb]

\includegraphics[width=0.49\textwidth]{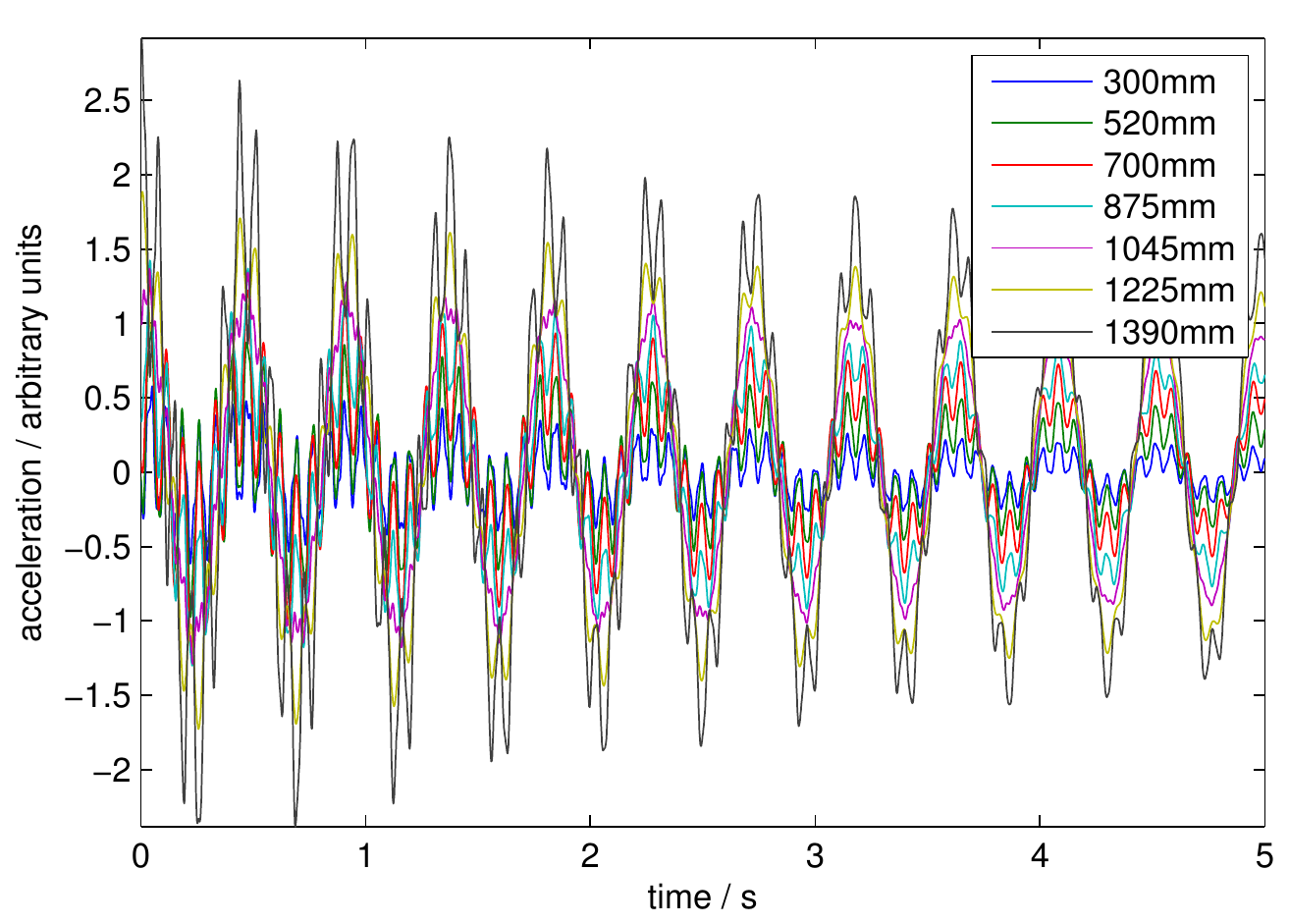}
\includegraphics[width=0.49\textwidth]{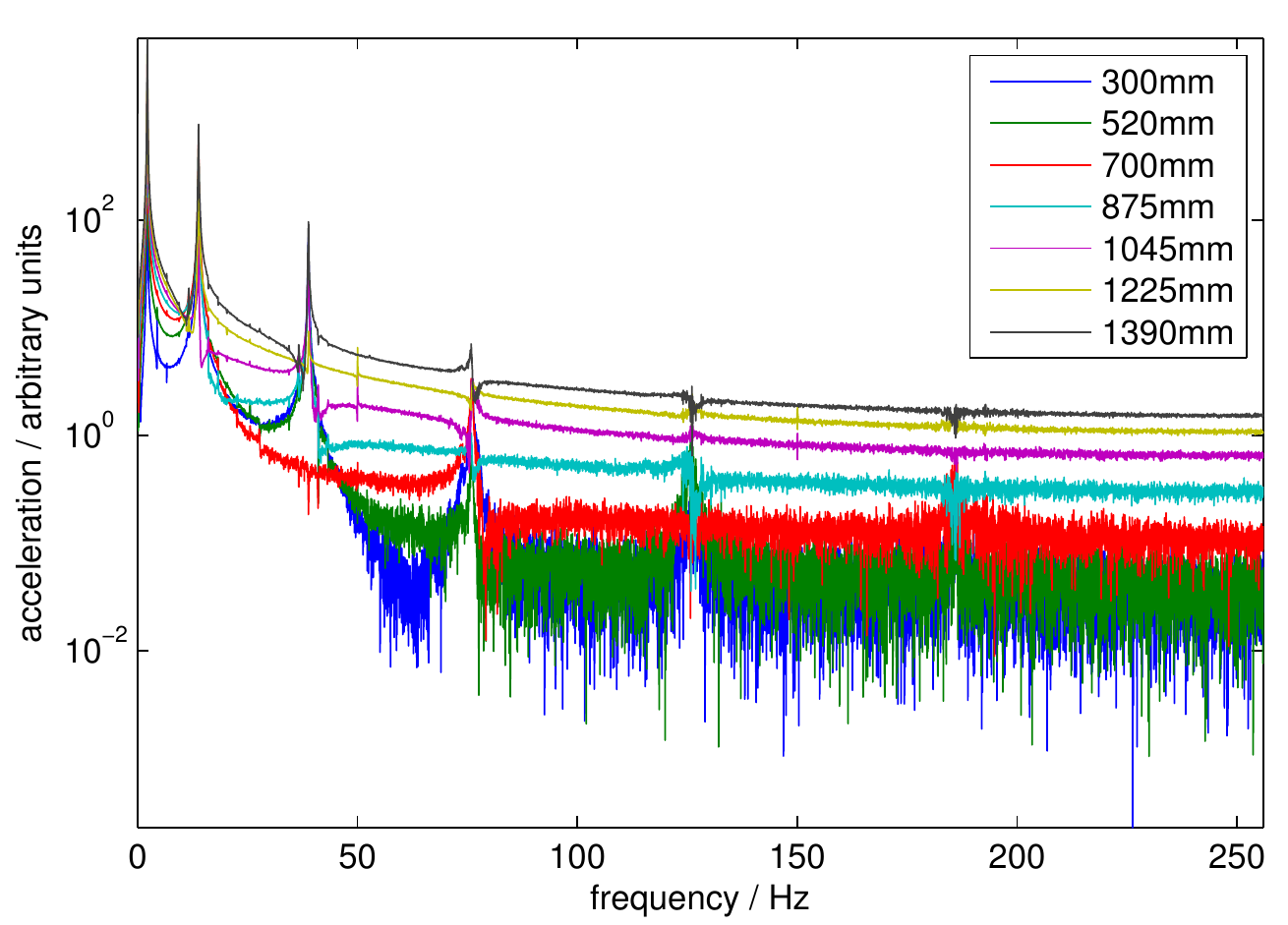}

\caption{Experimental data. Left: First 5 seconds of acceleration data. Right: Acceleration data in frequency domain.} \label{F:measured_data}

\end{figure}

Figure \ref{F:experiment} illustrates the measurement setup, including the location and orientation of the cut slot and the positions of the accelerometers. Figure~\ref{F:measured_data} shows a sample of the measured data in both the time  and frequency domains.

\subsection{Simulated measurement data} \label{S:simdata}

\begin{figure}[h]
\begin{center}
\includegraphics[scale=0.65]{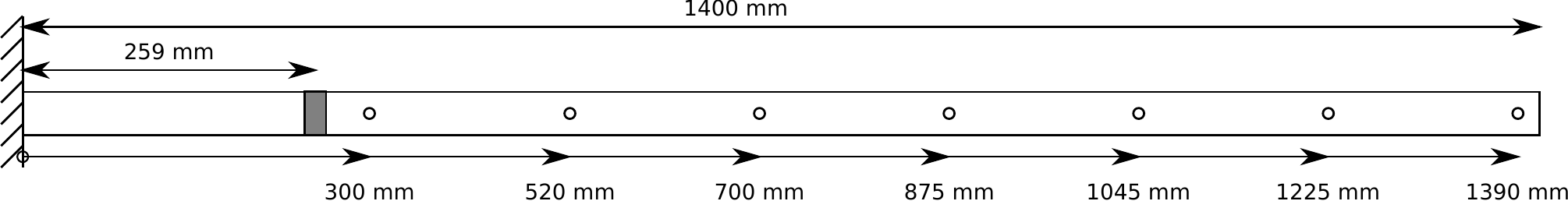}
\caption{Beam model used in data simulation. The small circles denote the points from which the simulated acceleration measurements were obtained. The gray rectangle represents element 19, which is the location of simulated damage.}
\label{F:simulation}
\end{center}
\end{figure}

The simulated data were setup to reflect the general conditions of the physical experiments including the location and severity of the damage. The sets of simulated data were generated in the time domain using the discretized model described in Section \ref{S:Model} using 100 elements of uniform size. The chosen nominal parameter values were $L=1400\ \text{mm}$, $EI = 131.25\ \text{N} \text{m}^2$, $\mu = 2.3\ \text{kg}/\text{m}$, $\alpha=0.15$ and $\beta=2 \cdot 10^{-5}$.
The simulated measurements were the accelerations measured at 7 points along the beam. Figure~\ref{F:simulation} illustrates the simulation model setup.

The initial conditions for the vibrations, modelling the bending of the beam by hand, were taken as the static solution of a displacement and a torque acting on the free end. The displacements were sampled from a $\mathcal{N}(10 \text{ mm},1 \text{ mm})$ distribution and the torques sampled from a $\mathcal{N}(0 \text{ Nm}, 0.5\text{ Nm})$ distribution. Environmental changes were simulated by randomly changing the length of the beam between each set of measurements according to temperature changes following a $\mathcal{N}(0\text{ K},5\text{ K})$ distribution. Each measurement set consisted of 30 seconds of data with a 512 Hz sampling rate of each of the seven measurement points. Normally distributed noise with standard deviation of 0.2\% of the initial signal amplitude was added in the time domain data to simulate a noisy measurement.

\begin{figure}[htbp]

\includegraphics[width=0.49\textwidth]{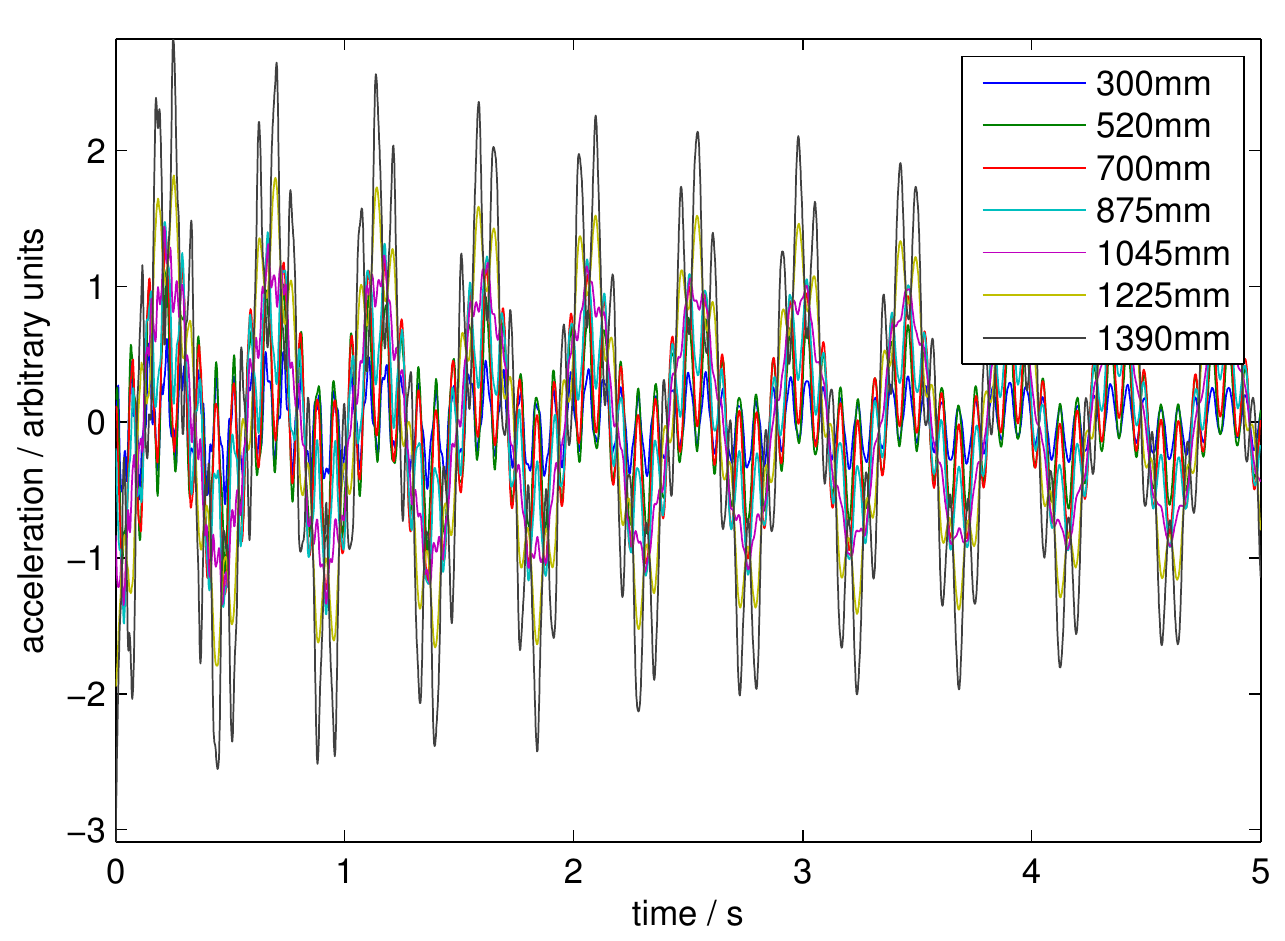}
\includegraphics[width=0.49\textwidth]{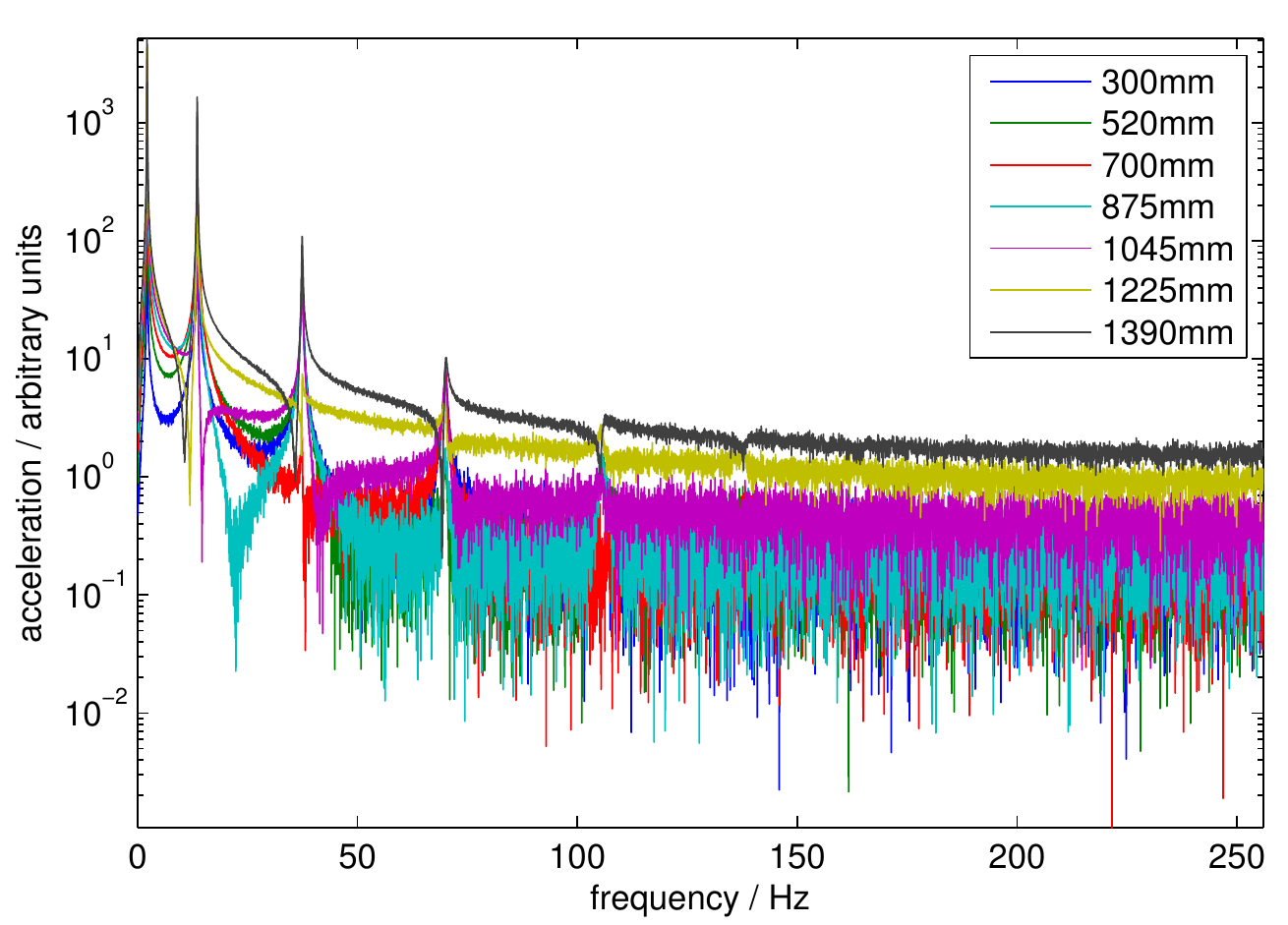}

\caption{Simulated data. Left: First 5 seconds of acceleration data. Right: Acceleration data in frequency domain.} \label{F:simulate_data}

\end{figure}

Five measurement sets were generated for the case where the beam is undamaged. Damage was then simulated using the damage model of Section~\ref{S:damage} with nonzero damage in element 19, which is centered around 259 mm from the fixed end of the beam. Four damage levels were tested, with the damage parameter value (in element 19) chosen as $0.125, 0.25, 0.375$ and $0.5$, and we refer to each of these damage levels as damage case 1, 2, 3, and 4, respectively. Five measurement sets were generated for each damage case for a total of 25 measurement sets including the undamaged case. Figure \ref{F:simulate_data} shows a sample of the simulated data in the time domain and in the frequency domain.% For each damage level, including no damage, 25 reconstructions of the damage parameters were obtained by both the EnKF and regularization schemes.

We emphasize that the data were simulated using a finer discretization of the beam model than is used by the UQ methods. We do this to try to avoid an inverse crime and to mimic the fact that an actual crack is generally narrower than the element size in any discretized model. As discussed in Section \ref{S:damage}, the severity of the damage is difficult to capture, but we are still generally able to identify with high probability the location of the most significant damage at or near the element containing the actual damage.

\subsection{Measurement pre-processing} \label{S:modal_identification}

The simulated and experimental data sets are in the form of acceleration data in the time domain. These data are used directly in the EnKF method, but the regularization and SACOM approaches are applied using modal data.
%In the case where experimental data is used within the framework of inverting observations of the computational model, we note that the measurement data comes from accelerometers in the form of time-signal data.

Significant pre-processing is required to obtain the vibration frequencies and observable mode shapes as described in Eq.~\eqref{eq:modal_solve}. One of the simpler ways to do this is to fit an autoregressive (AR) model to the data, and then solve for the poles of the model. With our data this approach gives a fitted model that tends to approximate the modes with high amplitude with multiple poles in their transfer function, while at the same time it tends to neglect the modes which have low amplitude. There are more elaborate pre-processing approaches, such as the Balanced Realization algorithm \cite{Benveniste1985}. Unfortunately, for our data it exhibits a similar problem as with the AR model. To overcome this, we exploit the knowledge that the beam is in free vibration, only lightly damped and that the modal frequencies are well separated as illustrated in Figures \ref{F:measured_data} and \ref{F:simulate_data}. As the beam is in free vibration, the model for each individual mode is a freely vibrating damped harmonic oscillator, given by
\begin{equation}
\xi_i(t) = A_i \exp(-t \zeta_i \omega_i) \cos(t \omega_i \sqrt{1-\zeta_i^2} + \phi_i),
\end{equation}
where $A_i, \omega_i, \zeta_i$ and $\phi_i$ are the vibration amplitude, vibration frequency, damping coefficient and phase shift of mode $i$, respectively.

Since each mode has a certain observable mode shape $\bm{h}_i$ and the modes are independent, the model seeks to explain the data $\bm{m}(t_k)$ as
\begin{equation}
\bm{y}(t_k) =  \sum_{i=1}^n \bm{h}_i \xi_i(t_k) \approx \bm{m}(t_k),
\end{equation}
where $n$ is the chosen number of modes.

The parameters $A_i, \omega_i, \zeta_i, \phi_i$ and the components of $\bm{h}_i$ are fitted from the measurement data through an optimization problem. Since the frequency response of each mode is band limited, we optimize in the frequency domain and concentrate only on the relevant frequency bins. More precisely,
\begin{equation}
(\bm{A}, \bm{\omega}, \bm{\zeta}, \bm{\phi}, \bm{h}) = \argmin_{(\bm{A}, \bm{\omega}, \bm{\zeta}, \bm{\phi}, \bm{h})} \sum_{j \in I} \| \bm{\hat{y}}(j) - \bm{\hat{m}}(j) \|^2.
\end{equation}
Here, $\bm{\hat{y}}(j)$ and $\bm{\hat{m}}(j)$ are the Fourier coefficients of the estimate and the data in frequency bin $j$, respectively. In the optimization, the observable mode shapes $\bm{h}_i$ are normalized so that $\bm{h}_i^T \bm{h_i} = 1$ and that the last component is non-negative. The set of relevant frequencies, denoted by $I$, are frequency bands which are centered around the peak frequencies in Figures \ref{F:measured_data} and \ref{F:simulate_data} and have an approximate -3 dB bandwidth. Concentrating on narrow frequency bands results in a significant reduction of noise in the data. The obtained pre-processed measurement is then collected as
\begin{equation}
\bm{m}_\omega = \begin{bmatrix}
\omega_1 \\
\bm{h}_1 \\
\vdots \\
\omega_n \\
\bm{h}_n
\end{bmatrix}
\end{equation}
matching the observation operator defined in Eq.~\eqref{eq:modal_solve}.

%We consider several measurements taken of each damage case, and run the pre-processing on each of the individual measurements. If the measurements and the pre-processing were ideal, the pre-processed measurements would be identical within each damage case. However, environmental effects as well as measurement noise cause deviations in the pre-processed measurements.

\subsection{Statistical parameters of the measurement noise} \label{S:measstatistics}

Both the regularization as well as the SACOM approach take measurement error into account through assuming a probability distribution for the measurement. This noise is assumed to be multivariate normal distributed, for which the mean and the covariance need to be estimated.

We consider several measurements taken of each damage case. If the measurements and the processing were ideal, the pre-processed measurements would be identical within each damage case. However, environmental effects, measurement error and model errors in the pre-processing cause deviations in the processed measurements. We take the measurement noise as additive and normal distributed, as was described in Eq. \eqref{eq:randmeasmodel}.
Assuming also that the measurement noise is independent of damage, we then see that the mean of the noise is given by
\begin{equation}
\bm{\mu} = E\left[\mathcal{E}\right] = E\left[\mathcal{M}|\mathcal{D}=\bm{d}\right]-f(\bm{d}).
\end{equation}
The only case where the damage state is really known is the undamaged case, for which $d=0$. The noise mean is thus estimated through the sample mean of the undamaged case measurements.
To estimate the noise coviariance, however, we use all of the measurement data.
\begin{equation}
\begin{aligned}
\bm{\Sigma} & = E\left[(\mathcal{E}-E\left[\mathcal{E}\right]) (\mathcal{E}-E\left[\mathcal{E}\right])^T\right]\\
& = E\left[ (\mathcal{M}-E\left[\mathcal{M}|\mathcal{D}=\bm{d}\right]) (\mathcal{M}-E\left[\mathcal{M}|\mathcal{D}=\bm{d}\right])^T|\mathcal{D}=\bm{d}\right].
\end{aligned}
\end{equation}
To compute this, we compute an estimate of $\bm{\Sigma}$ as the sample covariance of the measurements for each damage case individually. Notice, that this does not require knowledge of the value of $\bm{d}$, only that the value has remained constant for each sample. The final estimate of $\bm{\Sigma}$ is then obtained as the mean over all of the damage cases.

%In the measure-theoretic framework, the best results are obtained when the components of the data vector are {\em geometrically distinct} \cite{Butler2014a}. Conceptually, this means we prefer to use a data vector where the components define generalized contours that are most orthogonal to each other in the parameter space. For the sake of simplicity in presenting these results and to illustrate the benefits of working directly with the set-valued inverses, we use only the noisy data for the first frequency mode within this framework.

%The computational model is implemented as described in Section~\ref{S:fem} with $L=1.4$ m and $N+1=51$ uniformly spaced mesh points resulting in 50 uniformly sized elements $\set{T_i}_{i=1}^{50}$ of size $h=0.028$ m. The damage model is thus restricted to modeling a piecewise constant damage in terms of a stiffness reduction at a minimum length scale of $h$ as described in Section~\ref{S:damage}. In other words, the subspace of damage parameters used to form inputs to the computational model (see Figure~\ref{f:Abstract_flowchart}) is taken to be
%\begin{equation}\label{eq:d_space}
%    V_{h,d} := \set{\phi\in L^2(0,L)\, : \, \phi|_{T_i}\in \mathcal{P}_0(T_i) \text{ for } i=1,\ldots,50}.
%\end{equation}
%Here, $\mathcal{P}_0(T_i)$ denotes the space of constants on $T_i$. The time step for the implicit mid-point rule described in Section~\ref{S:forward1} is \textcolor{red}{Mika: What is it?}.

% ----------------------------------------------------------------
% ----------------------------------------------------------------
\section{Numerical results}\label{S:Numerics}
% ----------------------------------------------------------------
% ----------------------------------------------------------------

We considered two cases where the output data are obtained by either simulation of a refined numerical model or by experiments on a physically damaged beam. Each case is then sub-divided into the separate damage cases that are analyzed by each UQ method as discussed above.

The EnKF has the potential to run in real time, so it is used to provide an initial crude estimate of the damage and provide a numerical validation for the choice of penalty term in the regularization approach. The regularization results are subsequently used to define both a basis for a specified subdomain in the function space containing actual damage parameters and the domain of the {\em shape parameters} that parameterize this basis and define this subdomain in terms of the width, peak, and location of the damage. The measure-theoretic inverse is carried out with respect to these shape parameters to determine the probabilities of different types of localized damage. Since the damage is localized, we are particularly interested in quantifying uncertainties in the peak and locations as the regularization results suggest the width is bounded. Hence, we use the SACOM approach to produce marginal density plots of the joint peak and location parameters as well as marginal density plots of the width parameter for each damage case.

% ----------------------------------------------------------------
% ----------------------------------------------------------------
\subsection{Case I: Simulated observation data}
% ----------------------------------------------------------------
% ----------------------------------------------------------------

For the EnKF, we used the model described in Section~\ref{S:Model} discretized with 25 elements. The model dimensionality was intentionally kept low so that the EnKF can be used to quickly determine an initial estimate to the damage parameter. The model parameters not in the state vector were chosen as $L=1400 \text{ mm}$ and $\mu=2.3 \text{ kg}/\text{m}$. As the sampling frequency of 512 Hz is much faster than the dynamics of the beam, each pair of consecutive samples in the measurements are nearly identical. The assumption that EnKF relies on, that the measurement samples are statistically independent, thus does not hold. One solution would have been to ignore some of the measurement data, so as to have sufficient change between the samples, but as we do not have a very large amount of measurement data, we chose not to do this. Instead, we use the ensemble smoother approach (ES) \cite{EvensenBook} over 8 second windows in the data, and use all of it as the measurement. The size of the ensemble $N_e$ was 100 members. The initial ensemble was generated by first doing a least squares fit of the state vector to the first 8 seconds of the measurement data. Then each of the 100 members of the ensemble were generated by adding a random perturbation to the state vector.

For the regularization and SACOM methods, the beam model was discretized with 50 uniform elements. The model parameters were chosen as $L=1400 \text{ mm}, EI=133 \text{ N}\text{m}^2$ and $\mu = 2.3 \text{ kg}/\text{m}$. The measurements were pre-processed through the process of Section \ref{S:modal_identification} to obtain 3 of the lowest mode frequencies and observable mode shapes, for a total of 24 components in the measurement vector. The statistics of the measurement noise were estimated through the process of Section \ref{S:measstatistics}. The regularization method was found to be fairly robust with respect to the selection of the regularization parameter, giving meaningful results over a range of two orders of magnitude. The parameter value used in the presented results was chosen manually so that the obtained reconstructions looked smooth.

\subsubsection{No damage}

\begin{figure}[htb]
\begin{center}
\includegraphics[scale=0.6]{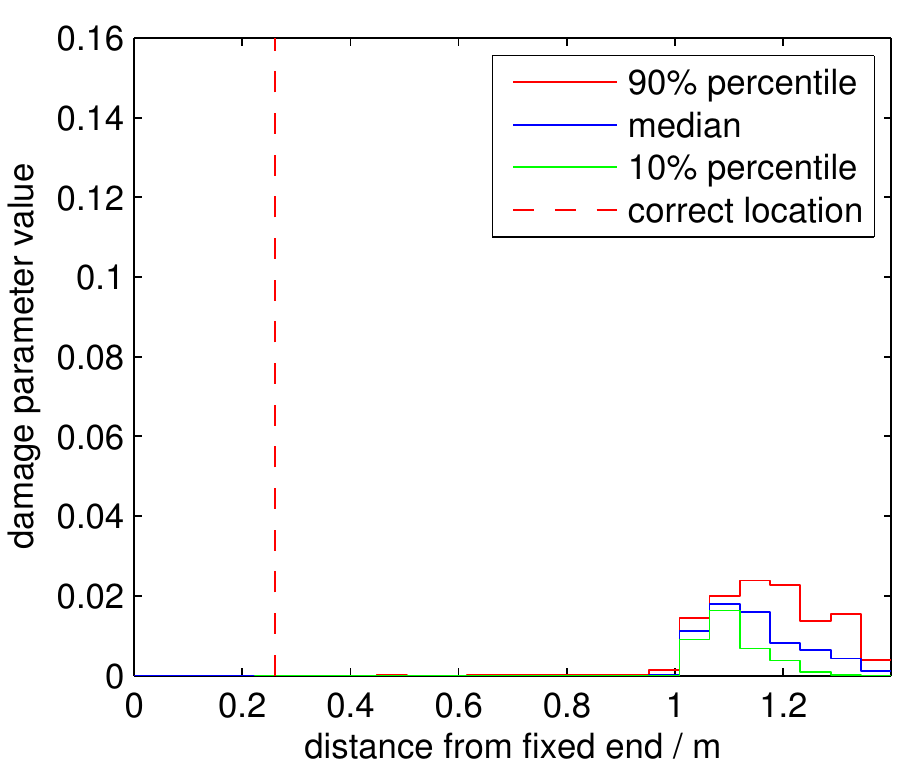}
\includegraphics[scale=0.6]{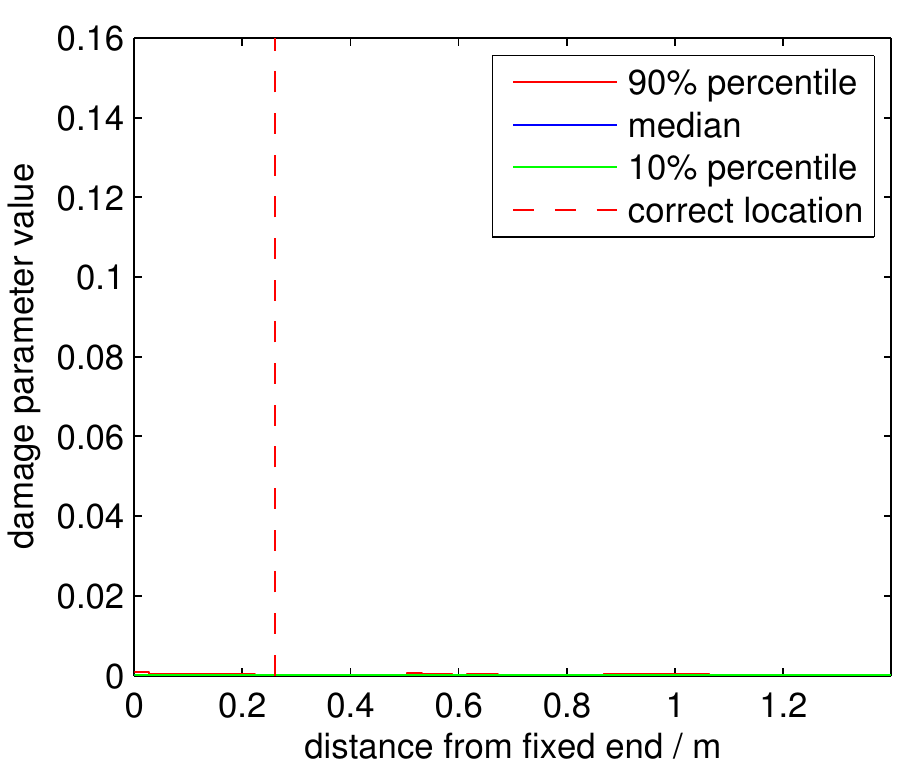}
\end{center}
\caption{Results for the simulated case with no damage. Left: statistics of the EnKF results. Right: statistics of the regularization method results.}\label{f:simData_d_case0}
\end{figure}

We run the EnKF and regularization methods to get a baseline of their performance when no damage is present. Damage parameter values were computed for each of the individual measurements using both methods. Figure~\ref{f:simData_d_case0} shows plots of damage parameters obtained using the EnKF and regularization methods. The plots show the statistics of each element-wise damage parameter independently over all of the obtained reconstructions, i.e. they give an idea of the ranges in which the parameter values in the reconstructions vary.

The EnKF results suggest there is a small amount of damage near the end of the beam. This indicates that the EnKF may have problems with detecting damage in this location when relatively small localized damages of the beam are present. The regularization method, however, identifies that any damage appears insignificant, which is a correct inference.

\subsubsection{Damage case 1}

\begin{figure}[htb]
\begin{center}
\includegraphics[scale=0.6]{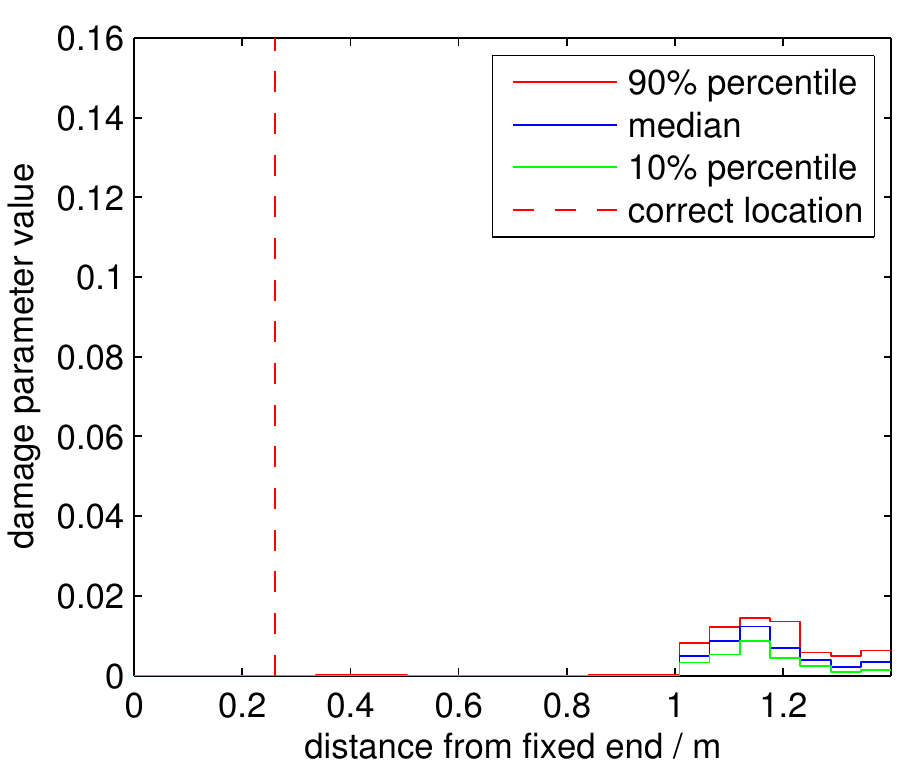}
\includegraphics[scale=0.6]{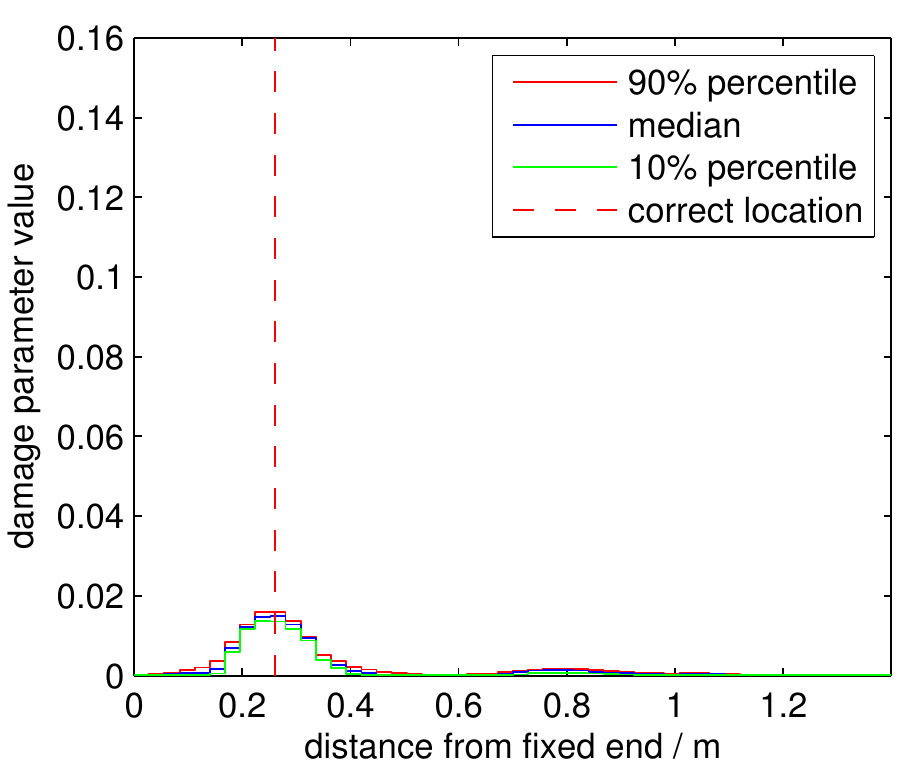} \\
\includegraphics[scale=0.6]{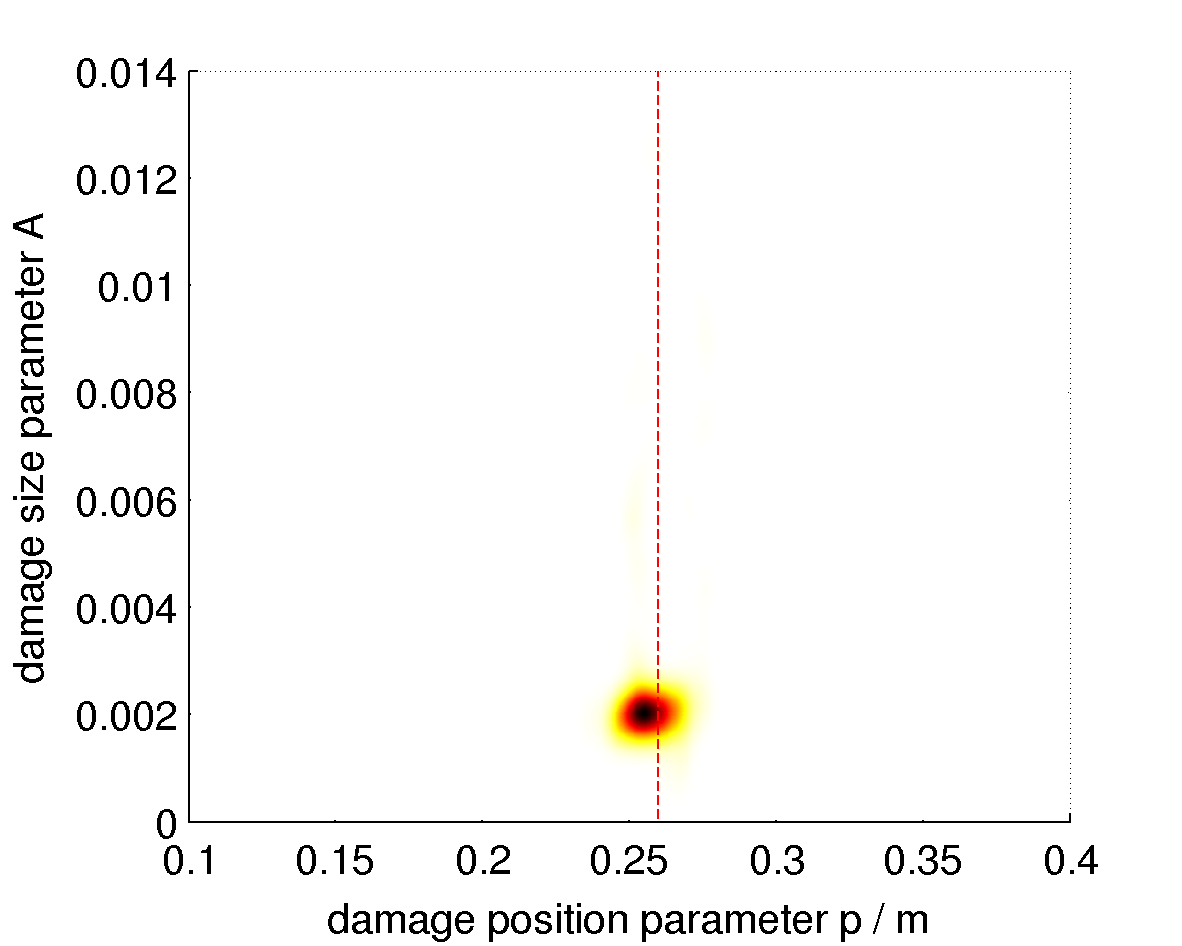}
\includegraphics[scale=0.6]{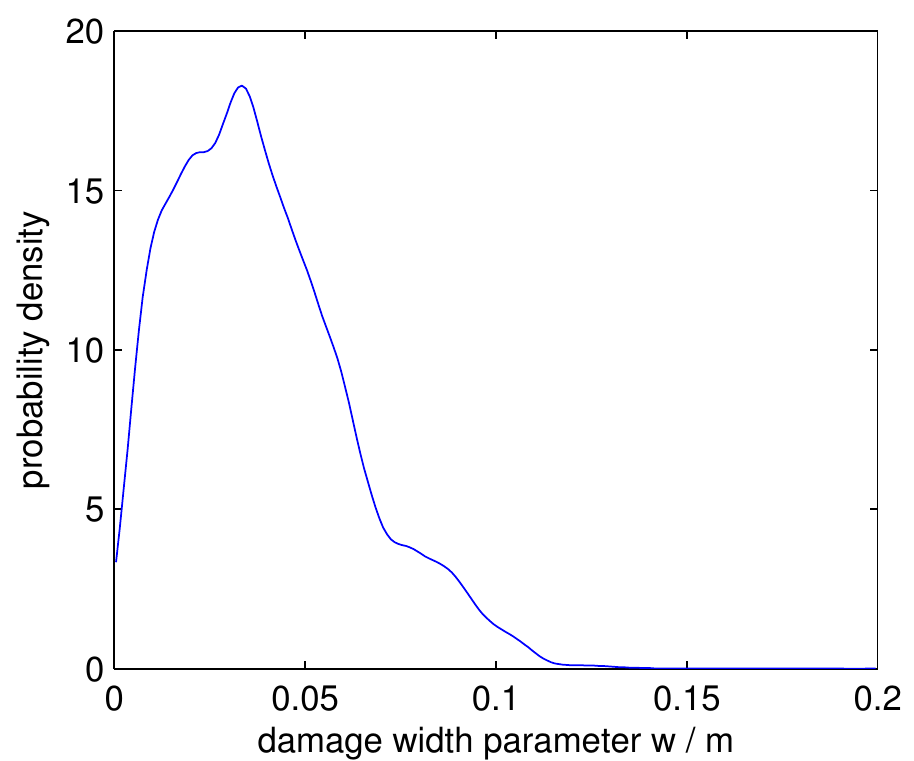}
\end{center}
\caption{Results for simulated damage case 1. Top left: statistics of the EnKF results. Top right: statistics of the regularization method results. Bottom left: marginal density of $(p,A)$ obtained with SACOM. The vertical red line represents the location of the simulated damage. Bottom right: marginal density of $w$ obtained with SACOM.}\label{f:simData_d_case1}
\end{figure}

In this case, the damage parameter of element 19 was set to value 0.125 in the simulation model. As before, damage parameter values were computed for each individual measurement using both the EnKF and the regularization methods. Figure~\ref{f:simData_d_case1} shows the statistics of results obtained from the EnKF and regularization, as well as SACOM, which we describe below.

The EnKF again suggests there is a small amount of damage near the end of the beam, which should be interpreted as a known incorrect detection, as the same output was observed with no damage. Moreover, no damage is detected at or near the correct location. However, the small amount of damage suggested by the EnKF can serve as a validation for the choice of penalty in the regularization and this is consistent among the remaining numerical results. The regularization method successfully detects the damage at the correct location. Compared with the undamaged case, we observe a significant change in the obtained parameter values.

As the regularization method results suggest that damage is present in a single location, we run the SACOM approach to further quantify the damage. Based on the regularization results, we take the damage to be a Gaussian function, which is characterized by three parameters: the overall damage size $A$, affected area width $w$ and the damage position $p$. The element-wise damage parameter vector $\bm{d}$ then takes the form
\begin{equation}
d_i = A \frac{1}{w \sqrt{\pi}} \exp\left[-(x_i-p)^2/w^2\right],
\end{equation}
where $x_i$ is the midpoint of the $i$th element. We can minimize the number of samples needed in SACOM by concentrating the parameters to the interesting range determined by the regularization results, i.e., we may restrict ranges of plausible values for $A$, $w$, and $p$ based on the regularization results. However, for the sake of simplicity in presenting consistent plots and exploring the probabilities of damage parameters at other physical locations, we allow $p$ to vary along the entire length of the beam, $w$ to vary between $[0,0.2]$ and $A$ to vary between $[0,0.015]$ which physically corresponds to no damage or the maximum observed peak of damage reconstructed by either the EnKF or regularization approaches for any case. We used the SACOM approach for each damage case (for both the simulated and experimental cases) with the same set of 5 million i.i.d.~uniform random samples in the Cartesian product of the specified intervals for $w$, $p$, and $A$ above. Since we are primarily interested in determining the location of damage, we restrict focus to the marginal density over the $(p,A)$ domain and to the marginal over $w$. We observe that the density for the position and damage size is very localized and that position is almost spot on with the correct location. The density function for the width is not as localized, but it is observed that it is unlikely for the damage to be very wide. These densities suggest that there is a small amount of damage in a fairly narrow localized region around the correct location of the simulated damage.

\subsubsection{Damage case 2}

\begin{figure}[htb]
\begin{center}
\includegraphics[scale=0.6]{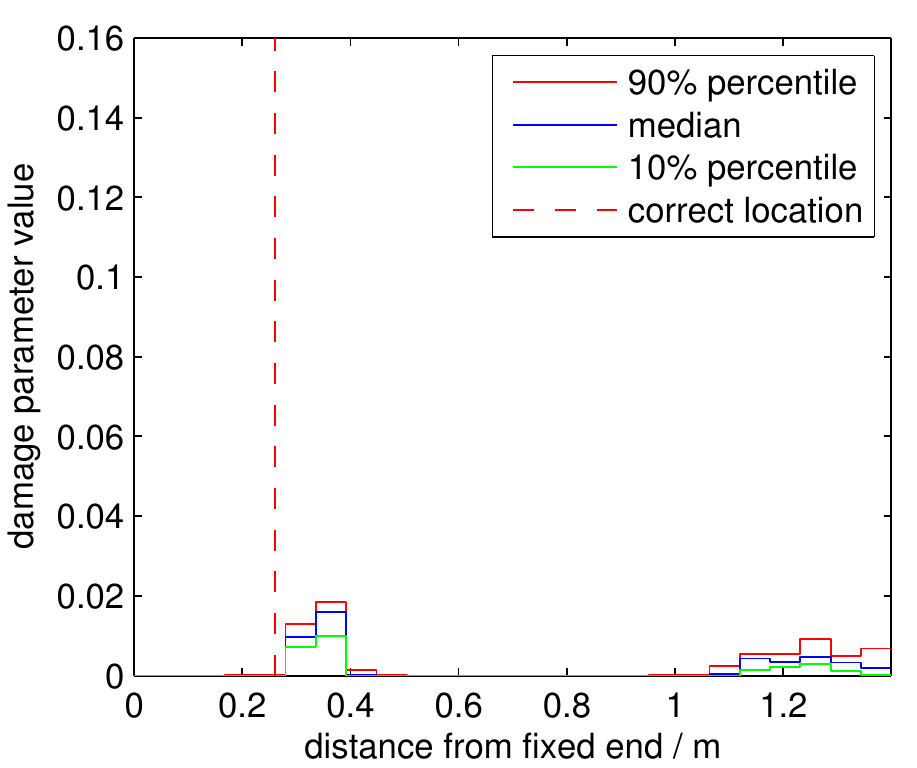}
\includegraphics[scale=0.6]{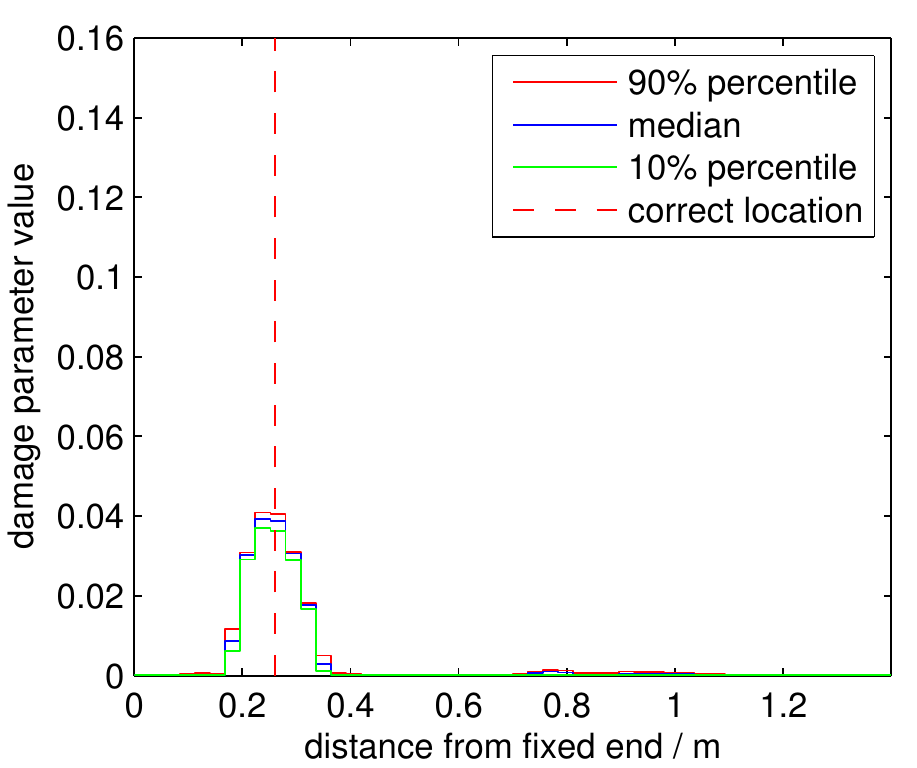} \\
\includegraphics[scale=0.6]{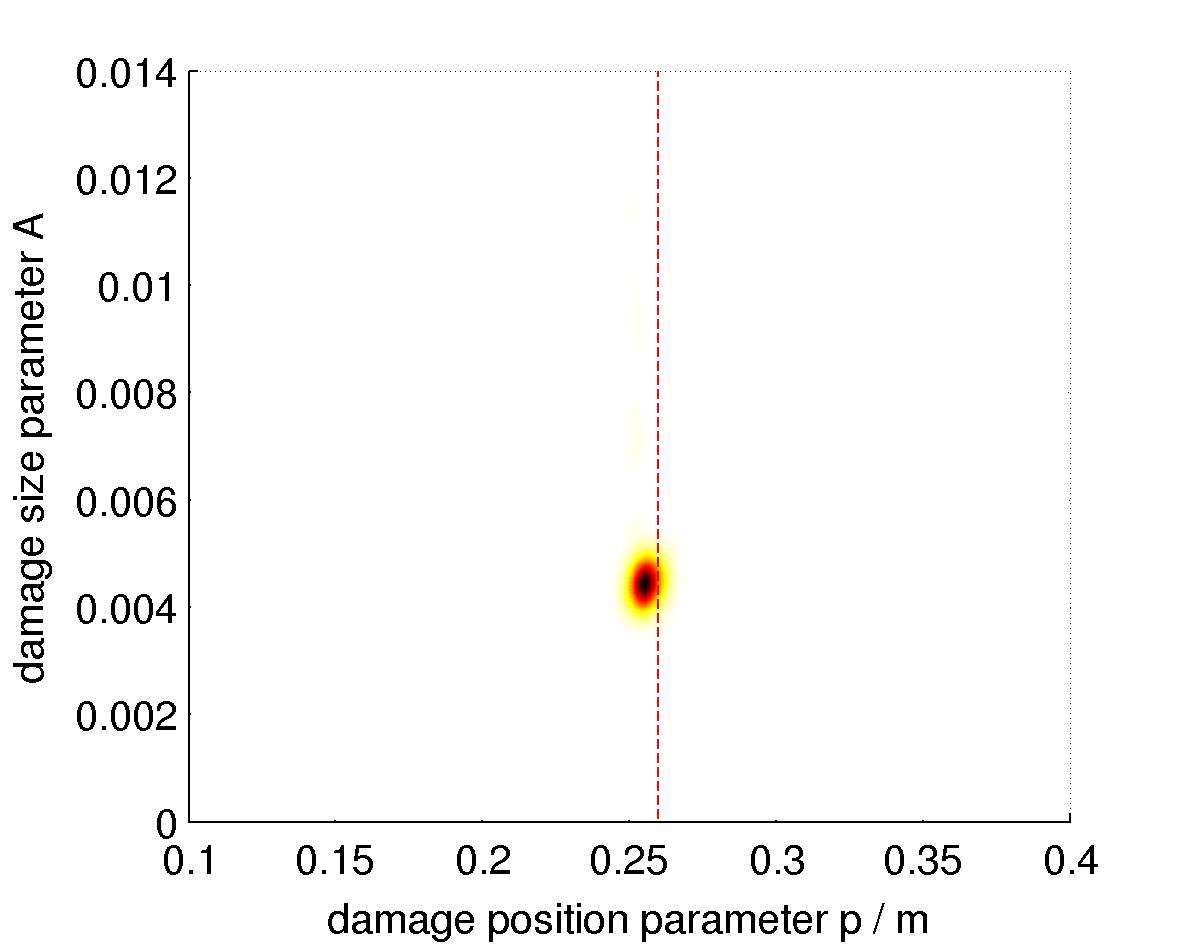}
\includegraphics[scale=0.6]{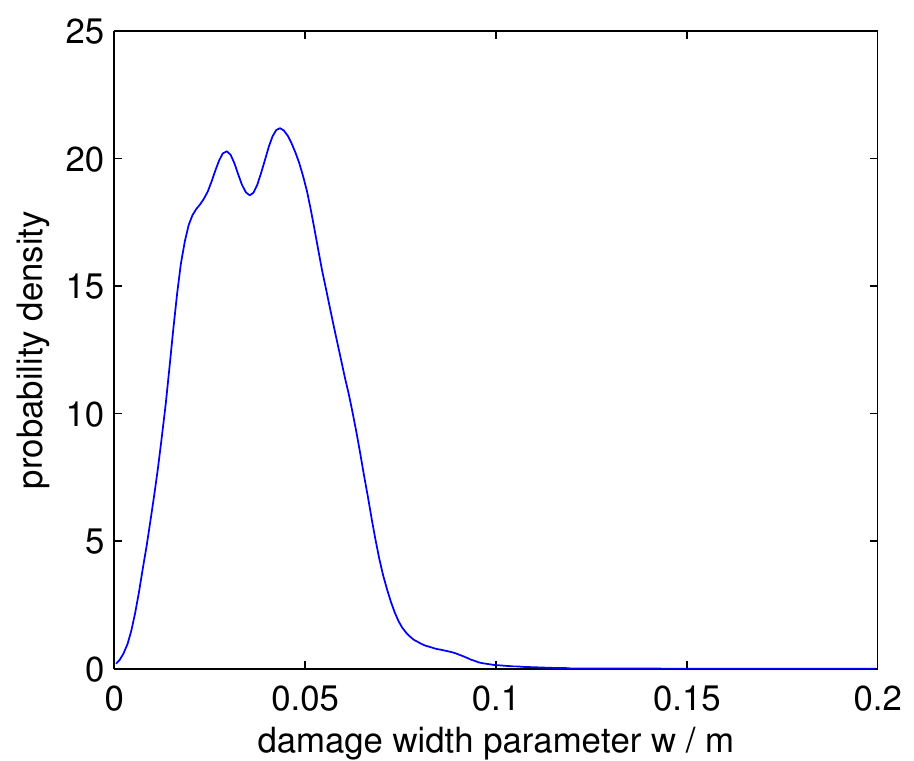}
\end{center}
\caption{Results for simulated damage case 2. Top left: statistics of the EnKF results. Top right: statistics of the regularization method results. Bottom left: marginal density of $(p,A)$ obtained with SACOM. The vertical red line represents the location of the simulated damage. Bottom right: marginal density of $w$ obtained with SACOM.}\label{f:simData_d_case2}
\end{figure}

In this case the damage parameter of element 19 was set to value 0.25 in the simulation model. Figure~\ref{f:simData_d_case2} shows the results of the various UQ methods for this damage case.

The EnKF continues to indicate a small amount of damage is present near the free end, which again should be interpreted as a known incorrect detection. However, the EnKF results now also indicate the presence of damage near the correct location. Although the amount of damage is in the same order as the spurious detections, it appears in a position where no spurious detections have previously been made and thus might be considered a positive detection. The regularization method again successfully detects the damage at the correct location, and the size of the peak is now larger than in the previous damage case. To further quantify the damage, we again use the SACOM approach assuming a single damage location based on the regularization results. We again observe that the marginal density over $(p,A)$ is very localized, with the most probable position parameters representing very accurately the actual damage location. The marginal for $w$ is again not quite as localized, but it is observed that large widths are very unlikely. These densities suggest that with high probability the damage is localized in a narrow area around the correct position, and that it is also more severe than the previous case.

\subsubsection{Damage case 3}

\begin{figure}[htb]
\begin{center}
\includegraphics[scale=0.6]{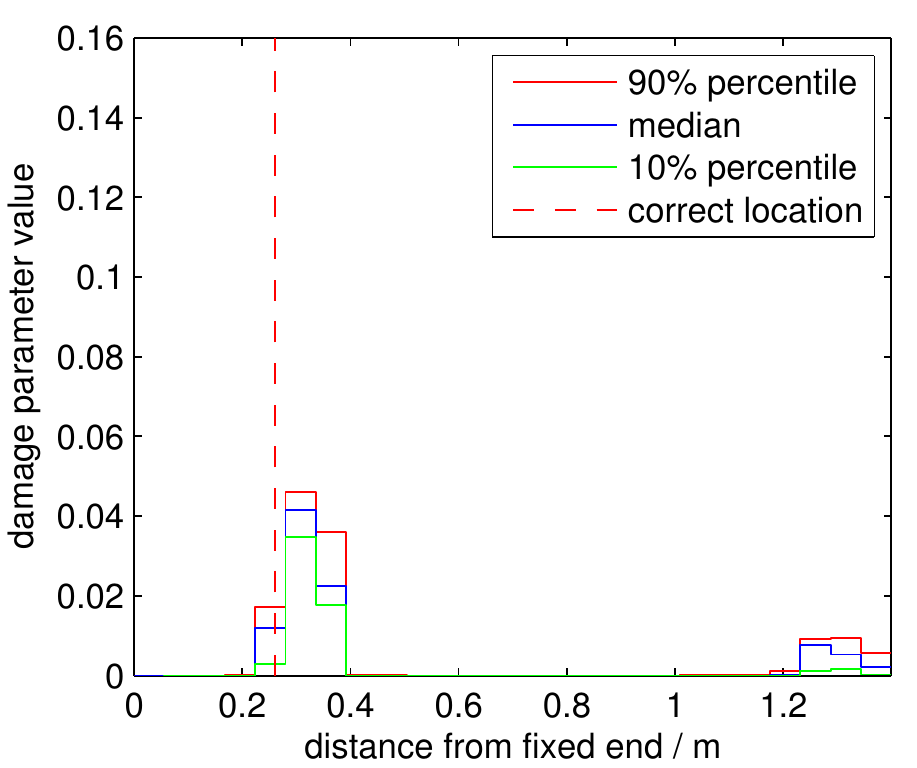}
\includegraphics[scale=0.6]{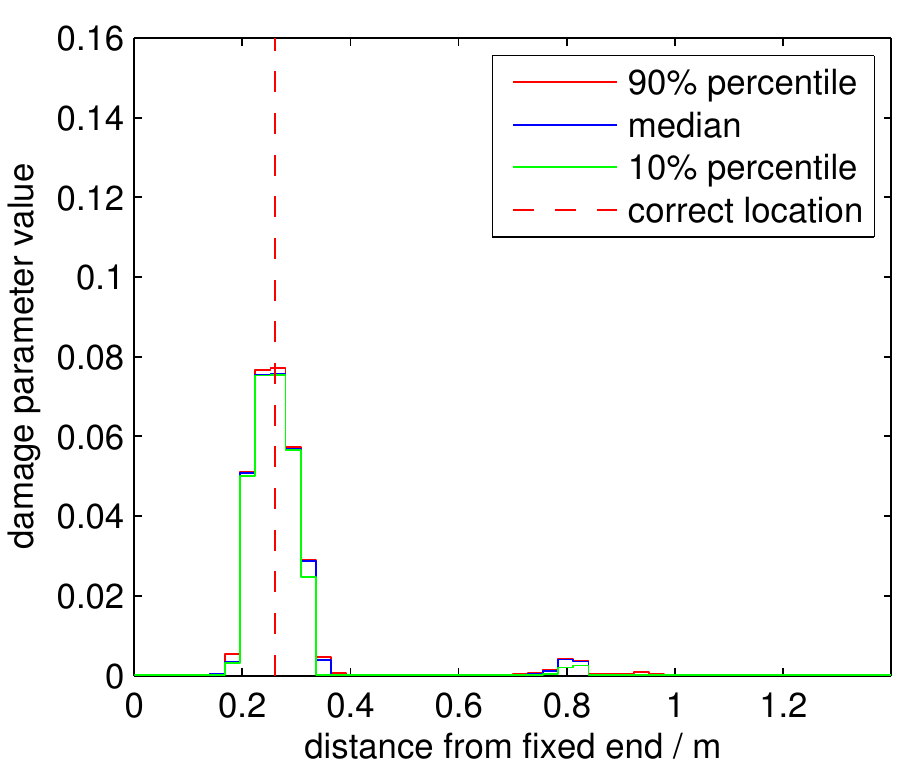} \\
\includegraphics[scale=0.6]{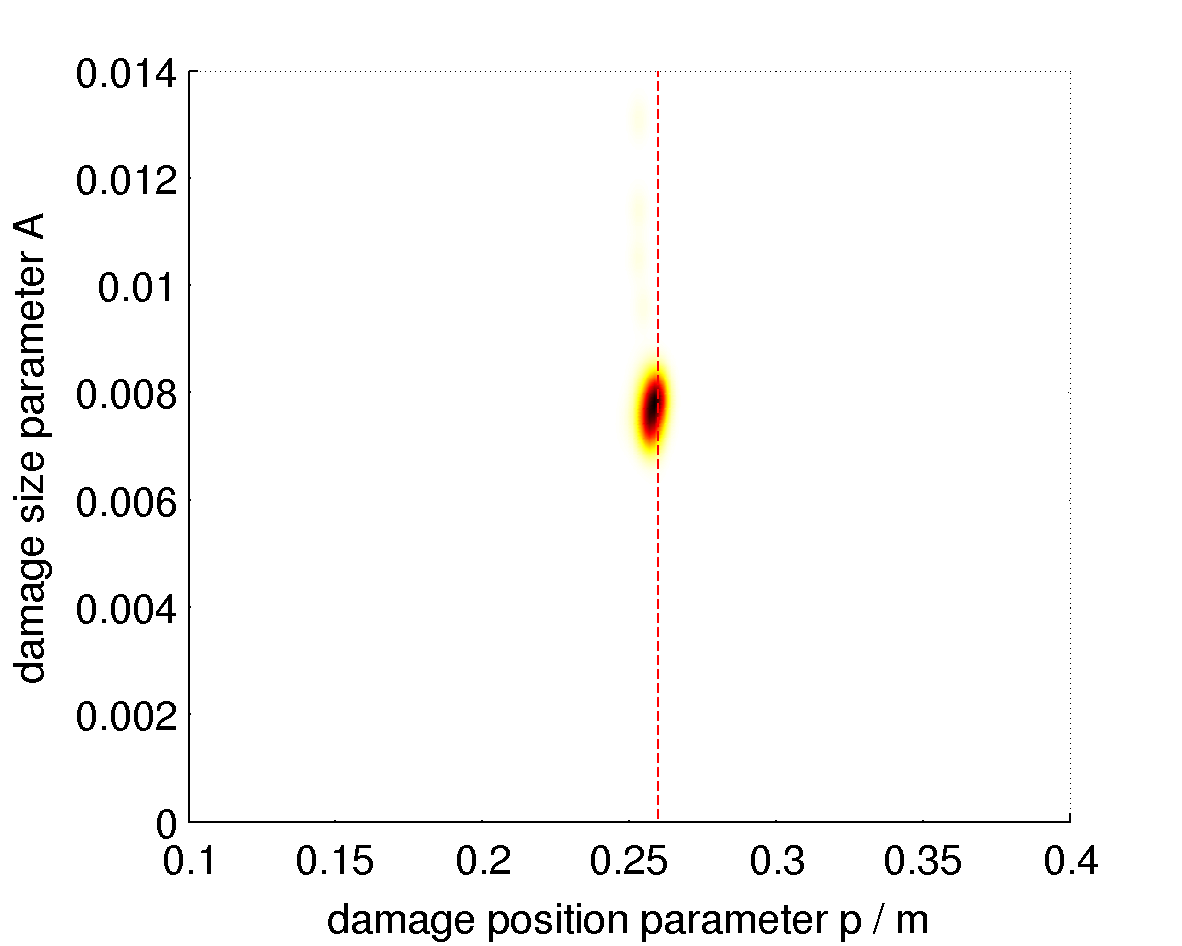}
\includegraphics[scale=0.6]{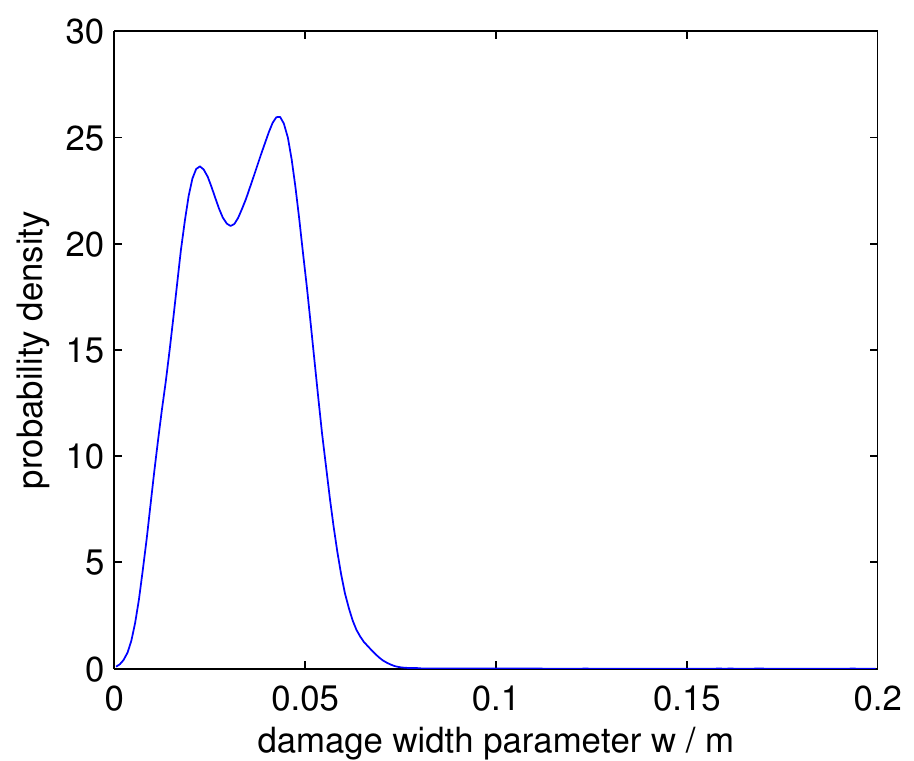}
\end{center}
\caption{Results for simulated damage case 3. Top left: statistics of the EnKF results. Top right: statistics of the regularization method results. Bottom left: marginal density of $(p,A)$ obtained with SACOM. The vertical red line represents the location of the simulated damage. Bottom right: marginal density of $w$ obtained with SACOM.}\label{f:simData_d_case3}
\end{figure}

In this case the damage parameter of element 19 was set to value 0.375 in the simulation model. Figure~\ref{f:simData_d_case3} shows the results of the various UQ methods for this damage case.

The EnKF now shows a clear indication of damage near the correct location. The damage size has increased from the previous case, and is now much larger than any of the spurious detections previously observed. The regularization method again detects damage at the correct location, and the amplitude has significantly increased from the previous damage case. Moreover, the SACOM results again suggest with high probability that there is an increase in the severity of the damage, and that it is localized around the true location with a fairly narrow width.

\subsubsection{Damage case 4}

\begin{figure}[htb]
\begin{center}
\includegraphics[scale=0.6]{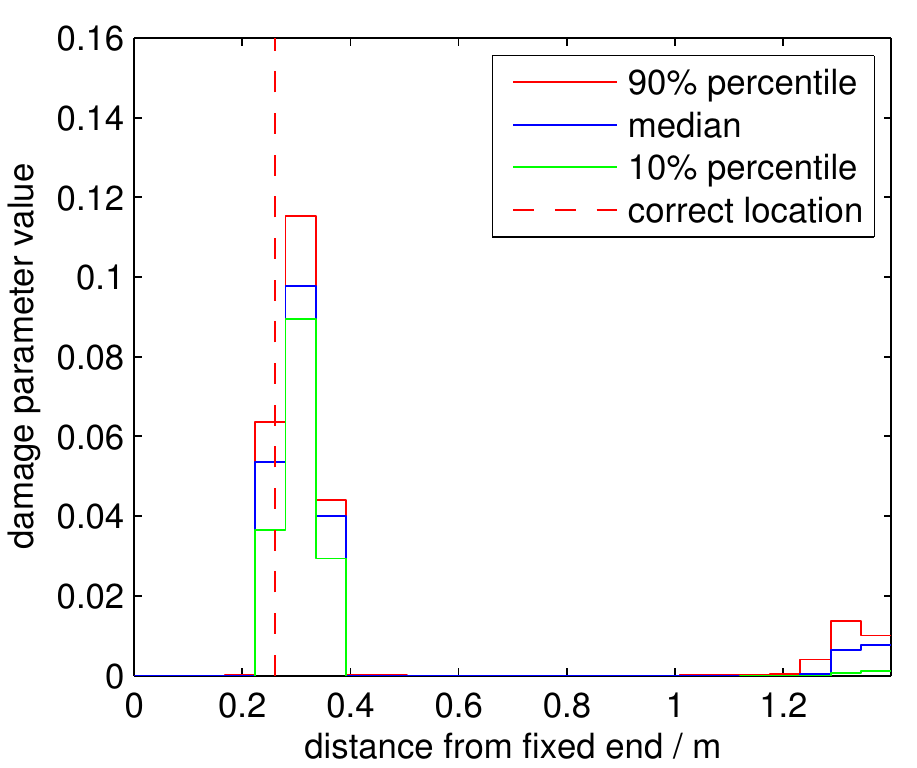}
\includegraphics[scale=0.6]{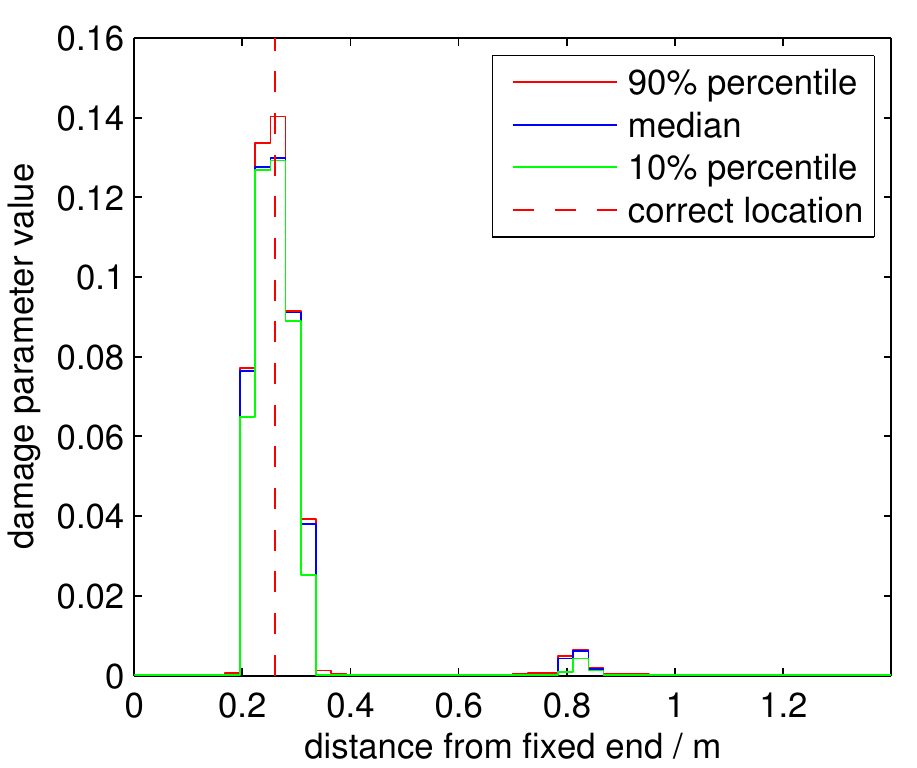} \\
\includegraphics[scale=0.6]{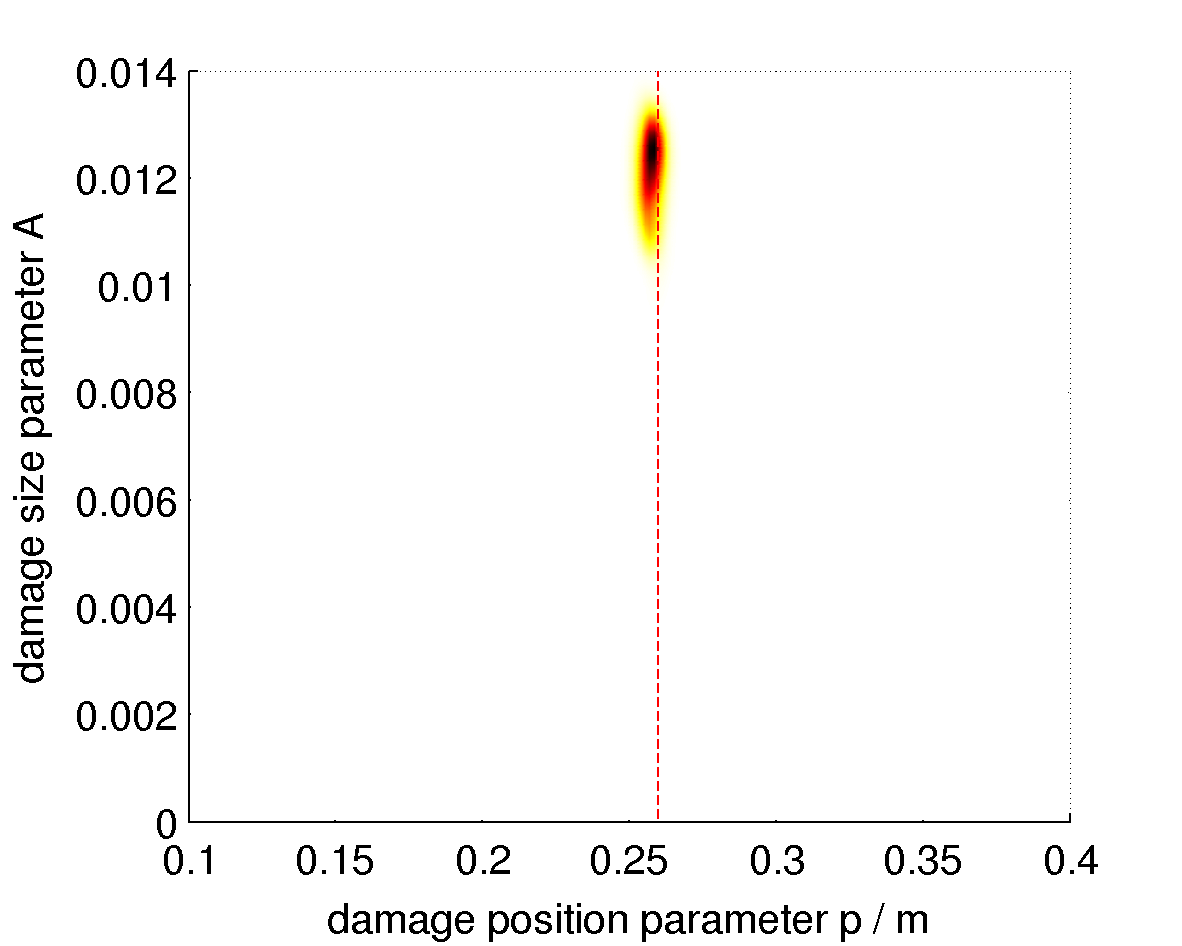}
\includegraphics[scale=0.6]{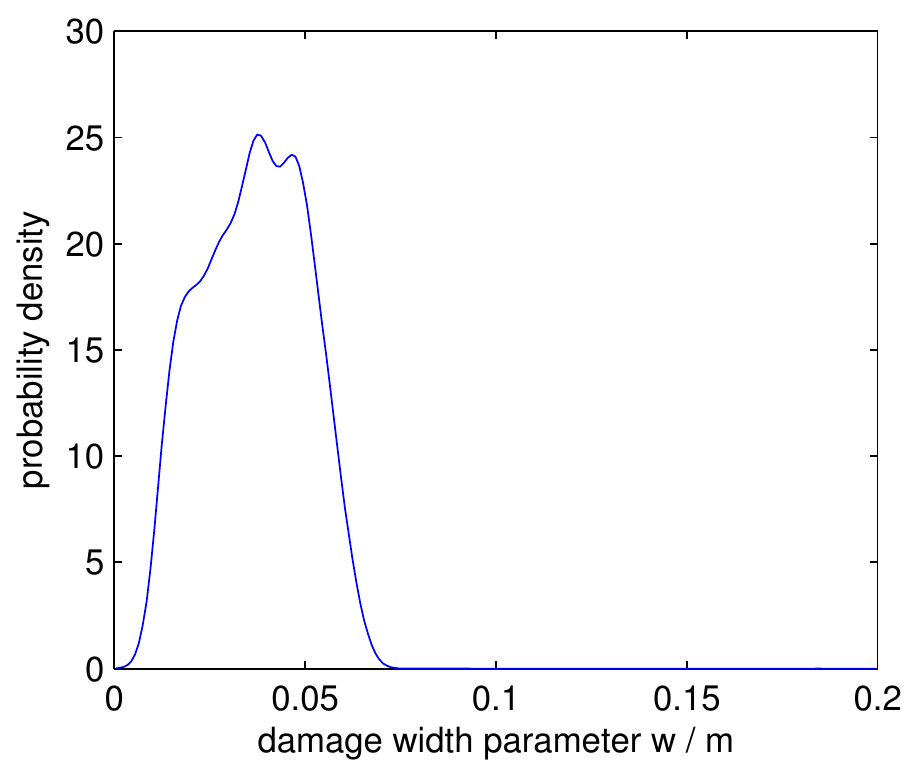}
\end{center}
\caption{Results for simulated damage case 4. Top left: statistics of the EnKF results. Top right: statistics of the regularization method results. Bottom left: marginal density of $(p,A)$ obtained with SACOM. The vertical red line represents the location of the simulated damage. Bottom right: marginal density of $w$ obtained with SACOM.}\label{f:simData_d_case4}
\end{figure}

In this case the damage parameter of element 19 was set to value 0.5 in the simulation model. Figure~\ref{f:simData_d_case4} shows the results of the various UQ methods for this damage case.

The EnKF again shows a clear indication of damage near the correct location, and the reconstructed damage size has increased from the previous damage case. The regularization method also continues to detect damage at the correct location, and with an amplitude that has also significantly increased from the previous case. The SACOM approach produces densities which suggests with high probability the damage is localized in a narrow area very near the correct location and that the damage has increased in amplitude from the previous cases.

% ----------------------------------------------------------------
% ----------------------------------------------------------------
\subsection{Case II: Experimentally obtained data}
% ----------------------------------------------------------------
% ----------------------------------------------------------------

As with the simulated data, the EnKF used the beam model with 25 elements, and the model parameter were chosen identically to the simulated data case. Again we used the ensemble smoother approach (ES) over 8 second windows as the measurement. The size of the ensemble $N_e$ was 100 members. The initial ensemble was generated by first doing a least squares fit of the state vector to the first 8 seconds of the measurement data. Then each of the 100 members of the ensemble were generated by adding a random perturbation to the state vector.

For the regularization and SACOM methods, the beam model was discretized with 50 uniform elements, and the model parameters were chosen identically to the simulated data case. The measurements were pre-processed to obtain 3 of the lowest mode frequencies and observable mode shapes and the statistics of the measurement noise were estimated. The regularization parameter value used in the presented results was chosen manually so that the obtained reconstructions looked smooth.

\subsubsection{No damage}

\begin{figure}[htb]
\begin{center}
\includegraphics[scale=0.6]{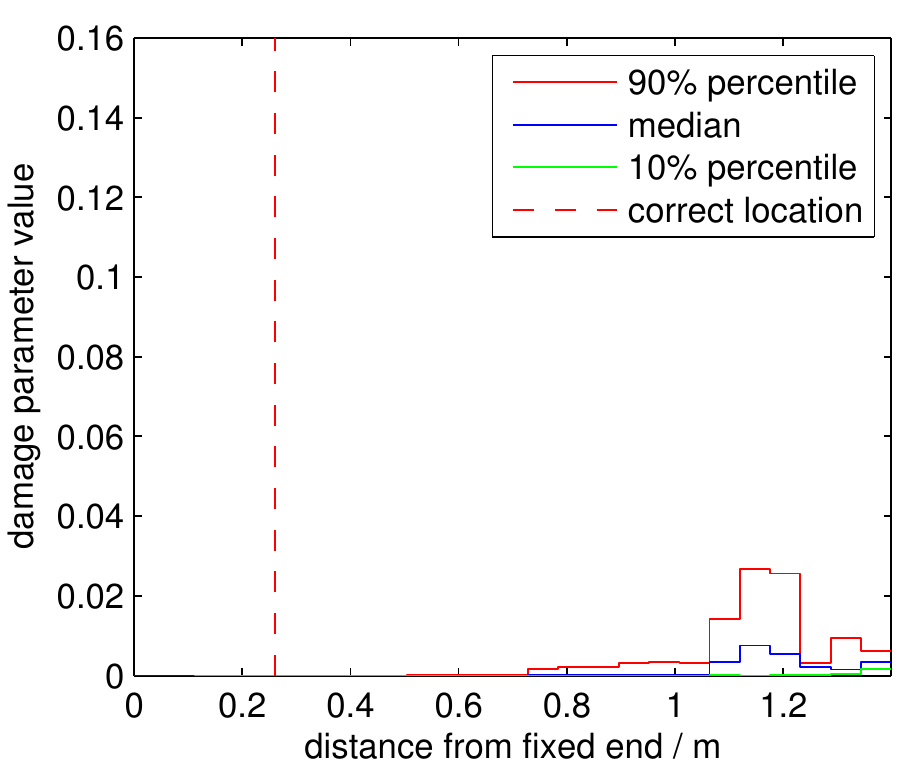}
\includegraphics[scale=0.6]{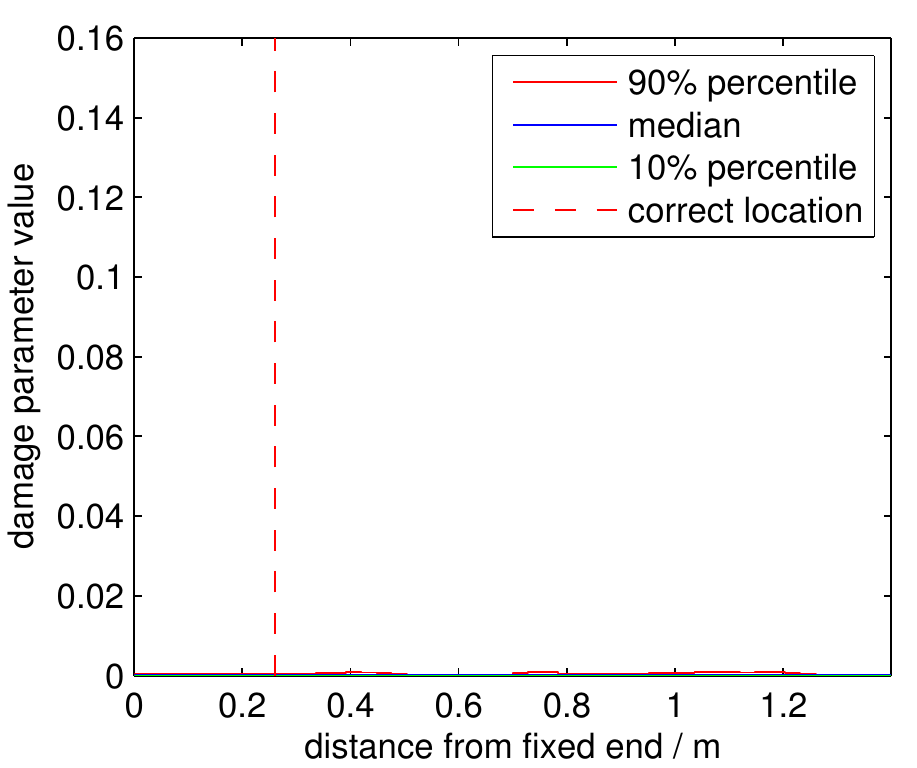}
\end{center}
\caption{Results for the experimental case with no damage. Left: statistics of the EnKF results. Right: statistics of the regularization method results.}\label{f:experData_d_case0}
\end{figure}

We again first test the EnKF and the regularization methods in the undamaged case to get a baseline of their performance. The damage parameter values are computed for each of the individual measurements. Figure \ref{f:experData_d_case0} shows statistics of the obtained results.

As with the simulated data, the EnKF incorrectly detects damage near the free end of the beam. The size of the spurious detection is approximately the same as for the undamaged case in the simulated data. The regularization method again correctly produces results with no significant damage. As the regularization method results suggest no damage is present, the SACOM approach is not used in this case.

\subsubsection{Damage case 1}

\begin{figure}[htb]
\begin{center}
\includegraphics[scale=0.6]{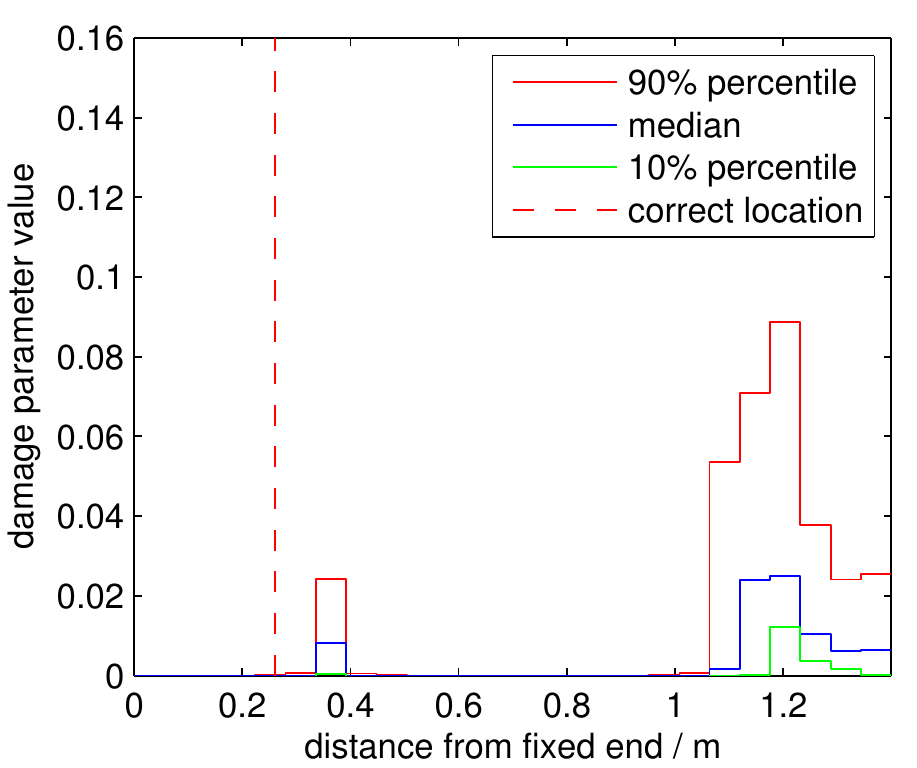}
\includegraphics[scale=0.6]{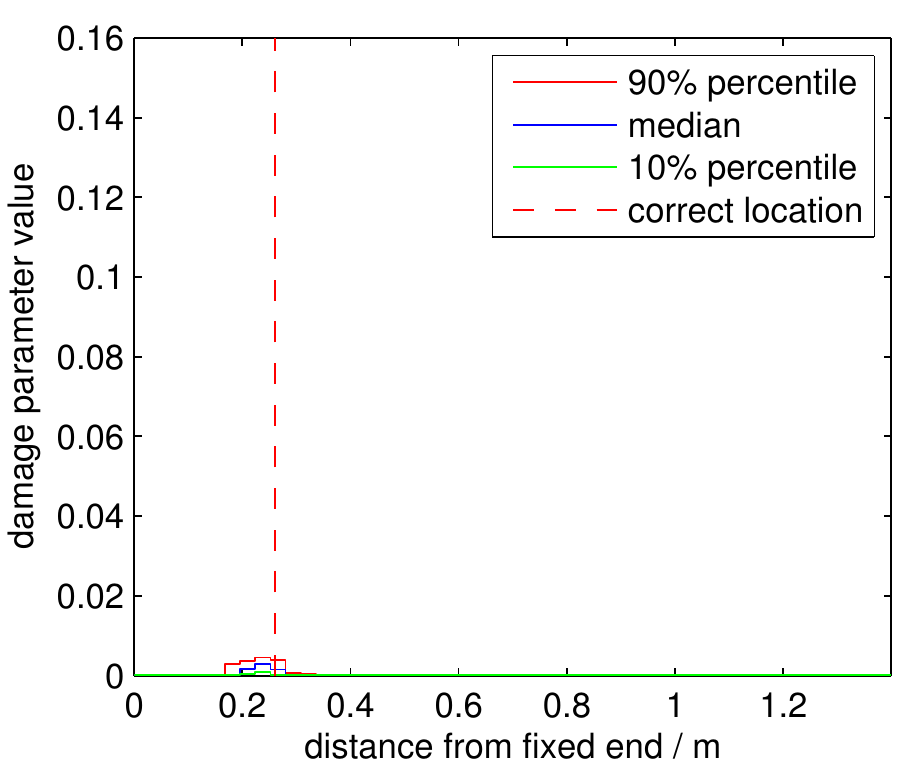} \\
\includegraphics[scale=0.6]{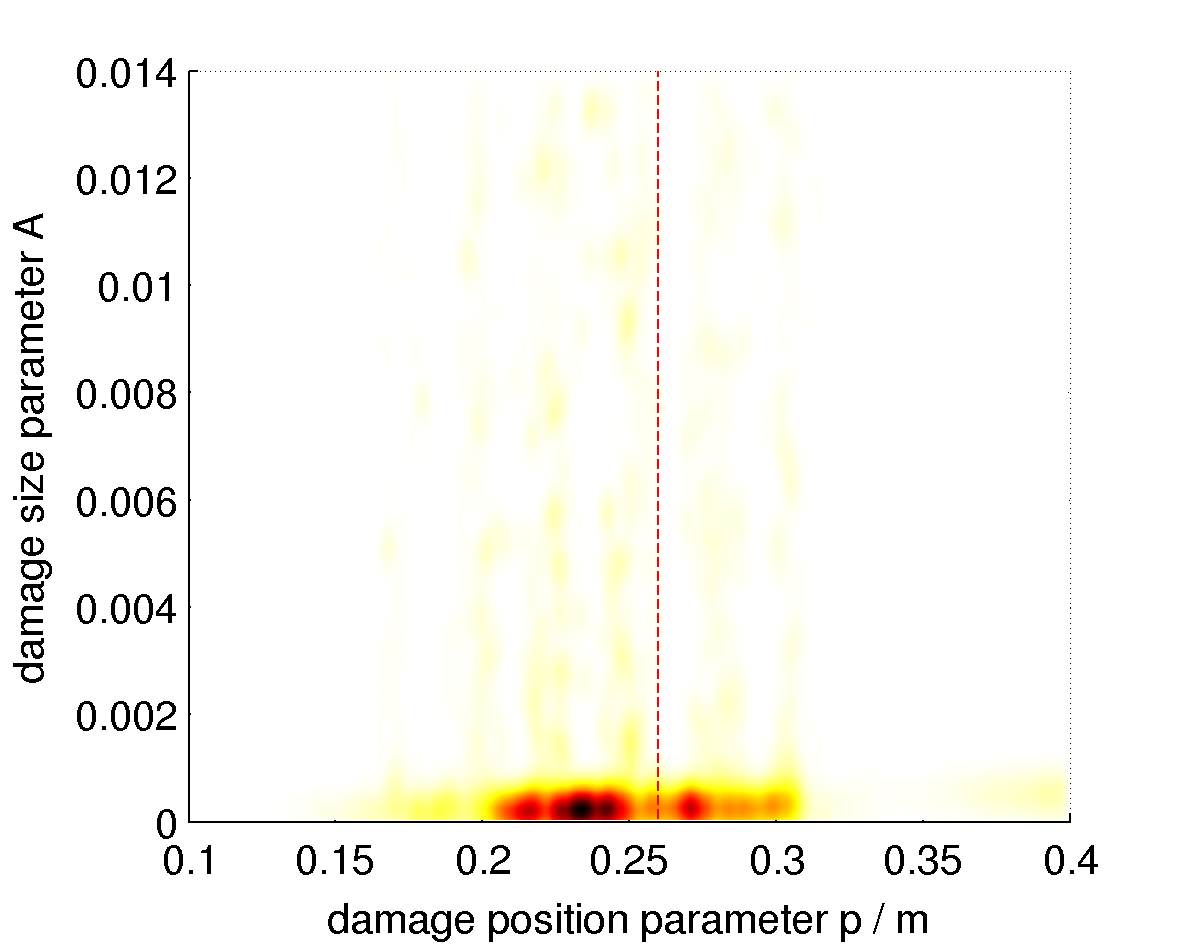}
\includegraphics[scale=0.6]{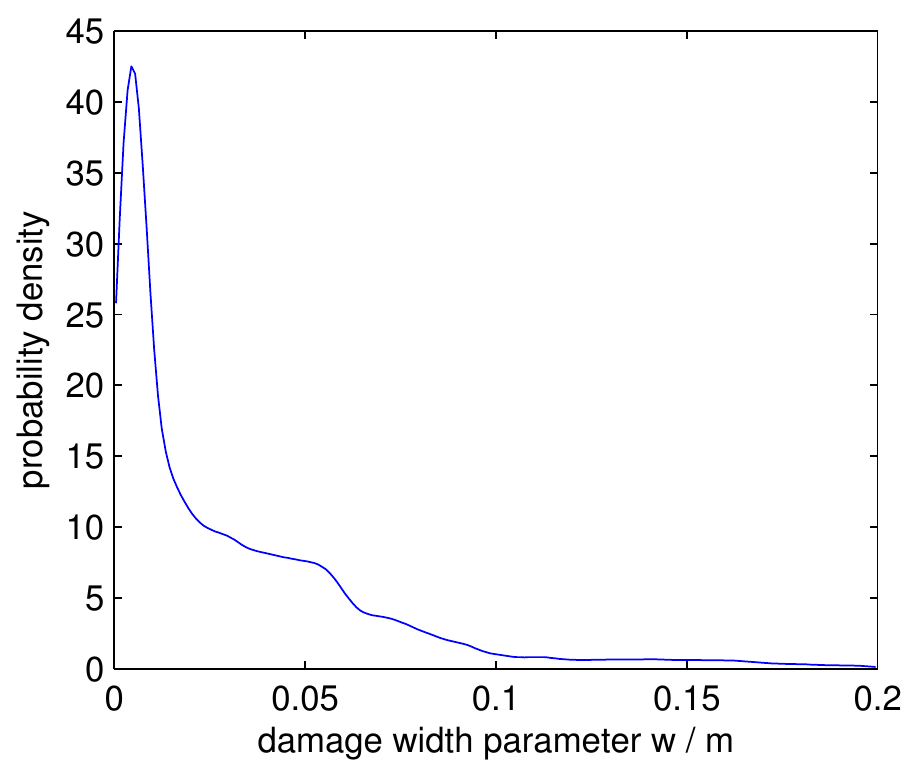}
\end{center}
\caption{Results for experimental damage case 1. Top left: statistics of the EnKF results. Top right: statistics of the regularization method results. Bottom left: marginal density of $(p,A)$ obtained with SACOM. The vertical red line represents the location of the afflicted damage. Bottom right: marginal density of $w$ obtained with SACOM.}\label{f:experData_d_case1}
\end{figure}

In this case the damage was a 5 mm deep and 1 mm wide slot at 260 mm from the fixed end of the beam. Damage parameters from the EnKF and regularization methods were computed for each of the individual measurements. Figure \ref{f:experData_d_case1} shows the statistics of these results, as well as the results obtained using the SACOM approach.

The EnKF again incorrectly detects damage near the free end of the beam. Due to the spurious damage detection in the undamaged beam, this should be interpreted as a known incorrect detection. A very small amount of damage is also detected near the correct location. However, the regularization method suggests a small amount of damage is present at the correct location. We use SACOM to further quantify the damage, which is observed to occur in a single location based on the regularization results. The marginal probability distribution for the damage position $p$ and damage size $A$ is not quite as localized as was seen with the simulated damage cases. While the mode of the distribution matches with the position obtained from the regularization method, SACOM gives a clear picture of the uncertainty involved. The damage size appears small, but it is greater than zero with high probability. The marginal density for the damage width $w$ implies that the affected area is very narrow. Together the densities suggest that with high probability there is a small amount of localized damage near the correct location.

\subsubsection{Damage case 2}

\begin{figure}[htb]
\begin{center}
\includegraphics[scale=0.6]{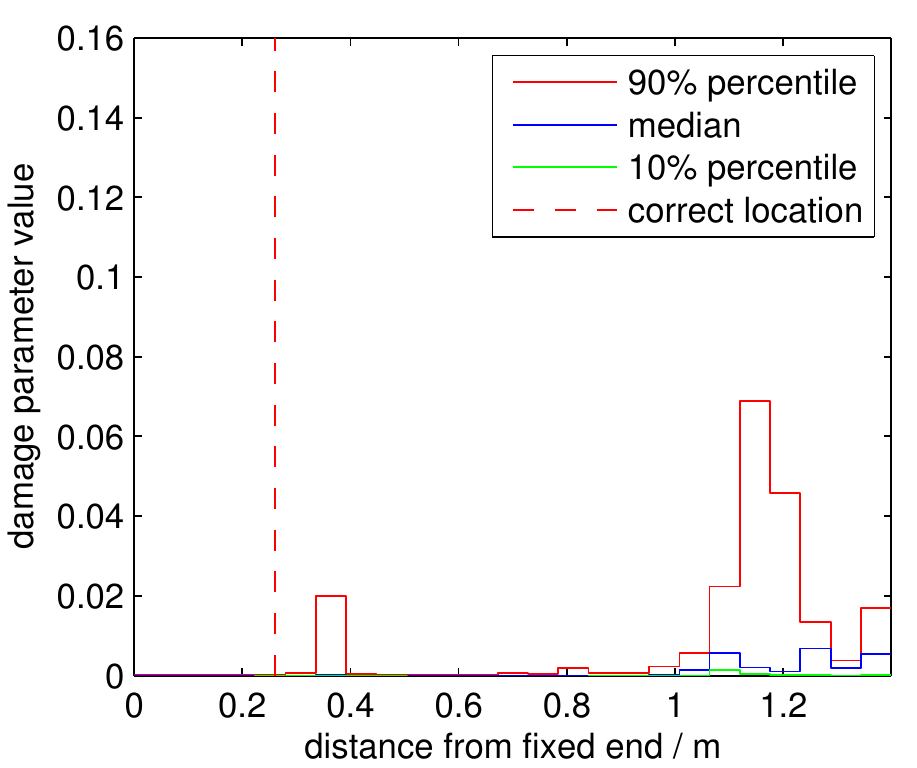}
\includegraphics[scale=0.6]{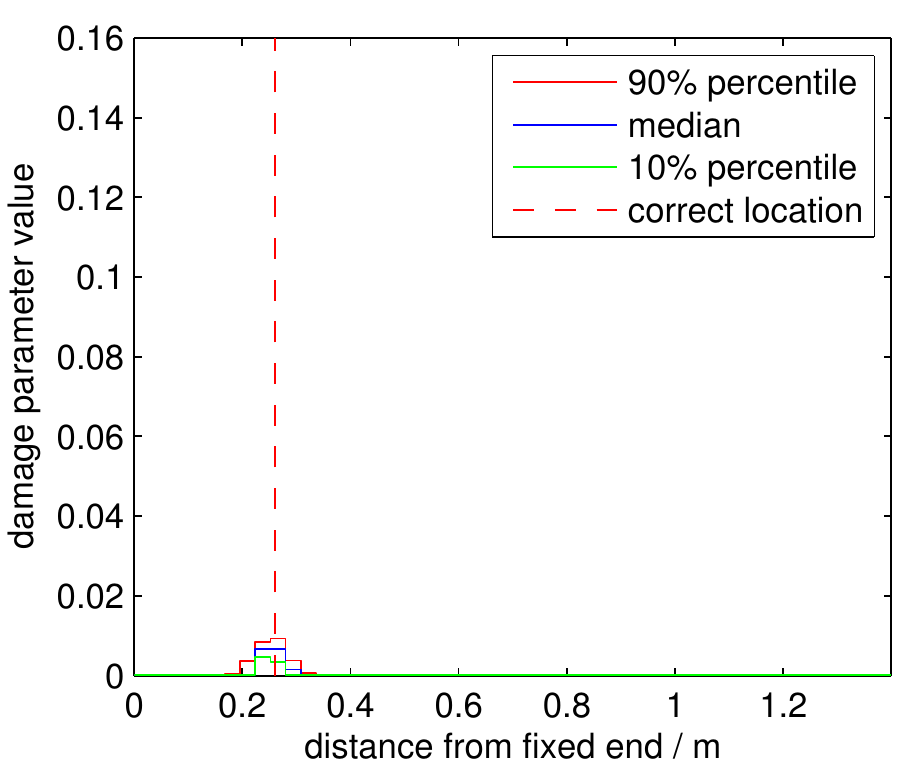} \\
\includegraphics[scale=0.6]{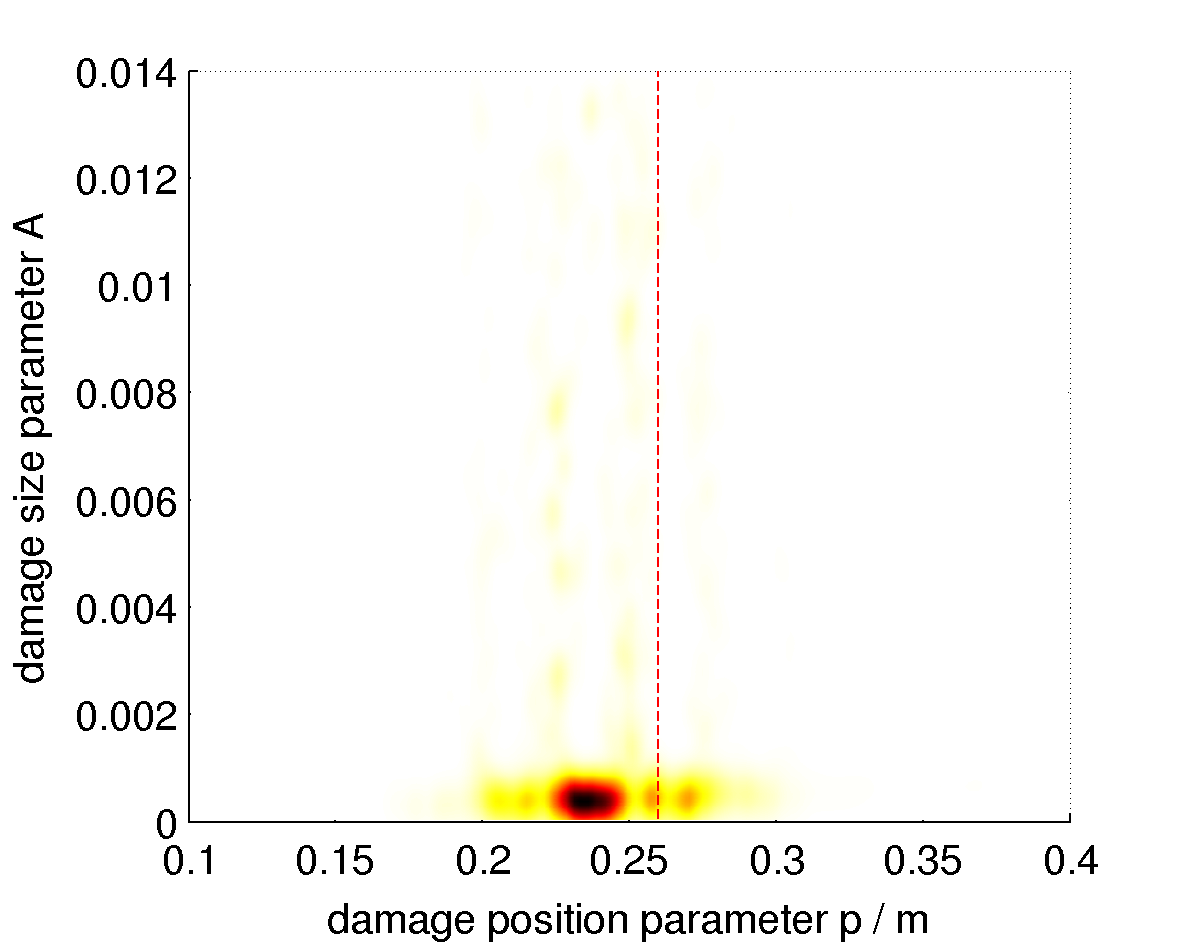}
\includegraphics[scale=0.6]{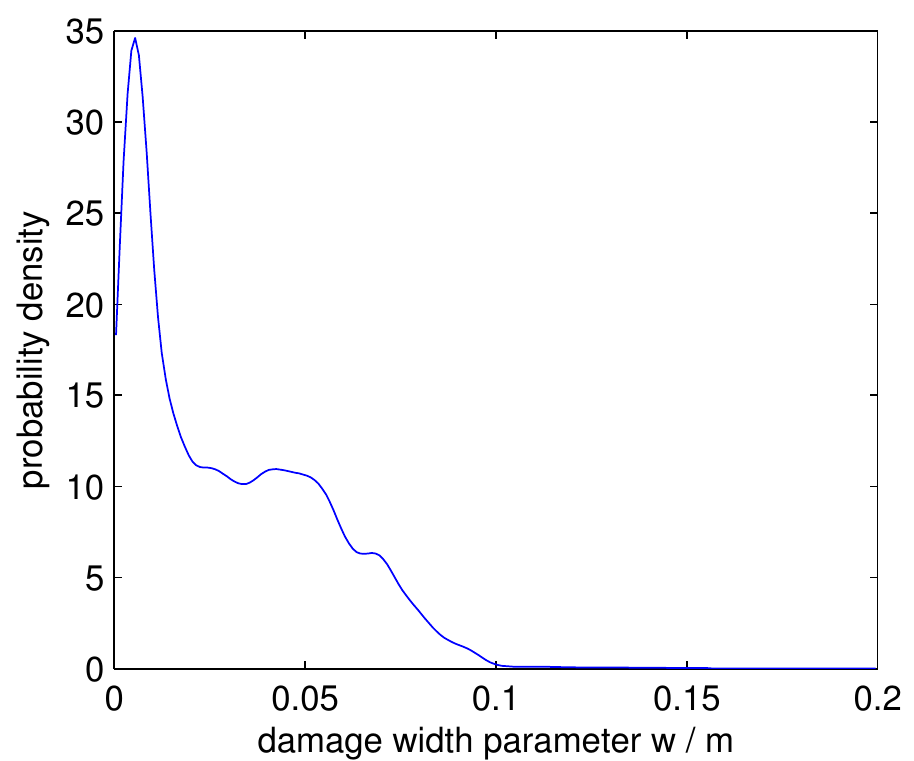}
\end{center}
\caption{Results for experimental damage case 2. Top left: statistics of the EnKF results. Top right: statistics of the regularization method results. Bottom left: marginal density of $(p,A)$ obtained with SACOM. The vertical red line represents the location of the afflicted damage. Bottom right: marginal density of $w$ obtained with SACOM.}\label{f:experData_d_case2}
\end{figure}

In this case the damage was a 10 mm deep and 1 mm wide slot at 260 mm from the fixed end of the beam. Figure \ref{f:experData_d_case2} shows the results from the three methods.

The EnKF incorrectly detects damage near the free end of the beam, which is again interpreted as a known incorrect detection. A very small amount of damage is also detected near the correct location. The regularization again detects a small amount of damage at the correct location, and the detected damage size is now larger than in the previous damage case. The SACOM approach is again used to further quantify the damage, which the regularization results suggest is in a single location. The marginal probability density of $(p,A)$ again suggests with high probability that the damage is located near the correct location. Moreover, the severity of the damage has increased compared to the previous case, and the affected area remains quite narrow with high probability.

\subsubsection{Damage case 3}

\begin{figure}[htb]
\begin{center}
\includegraphics[scale=0.6]{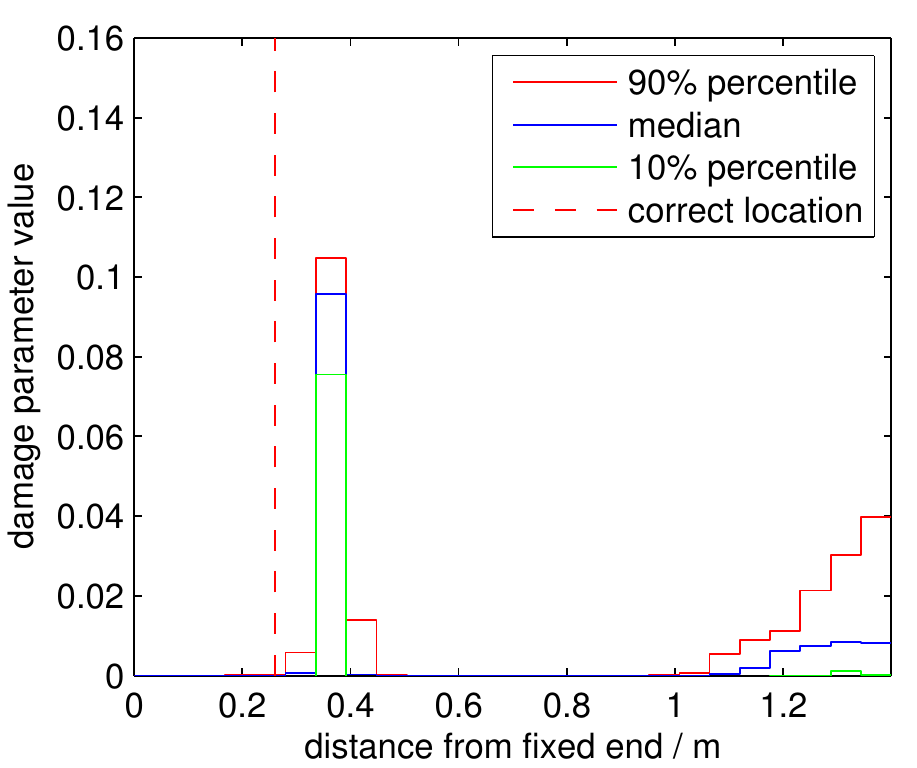}
\includegraphics[scale=0.6]{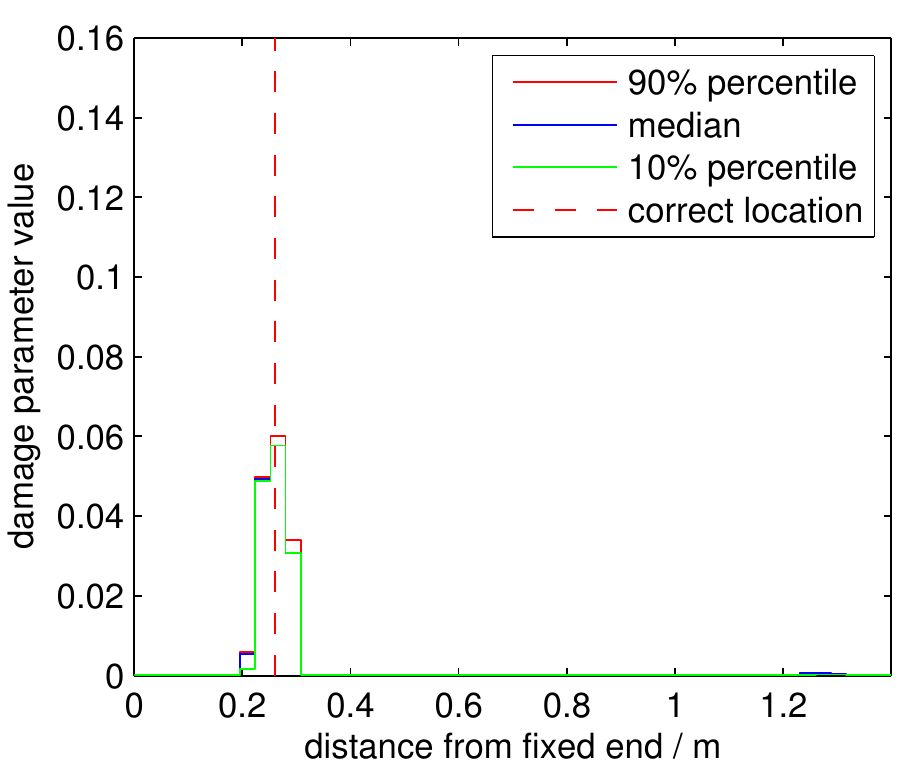} \\
\includegraphics[scale=0.6]{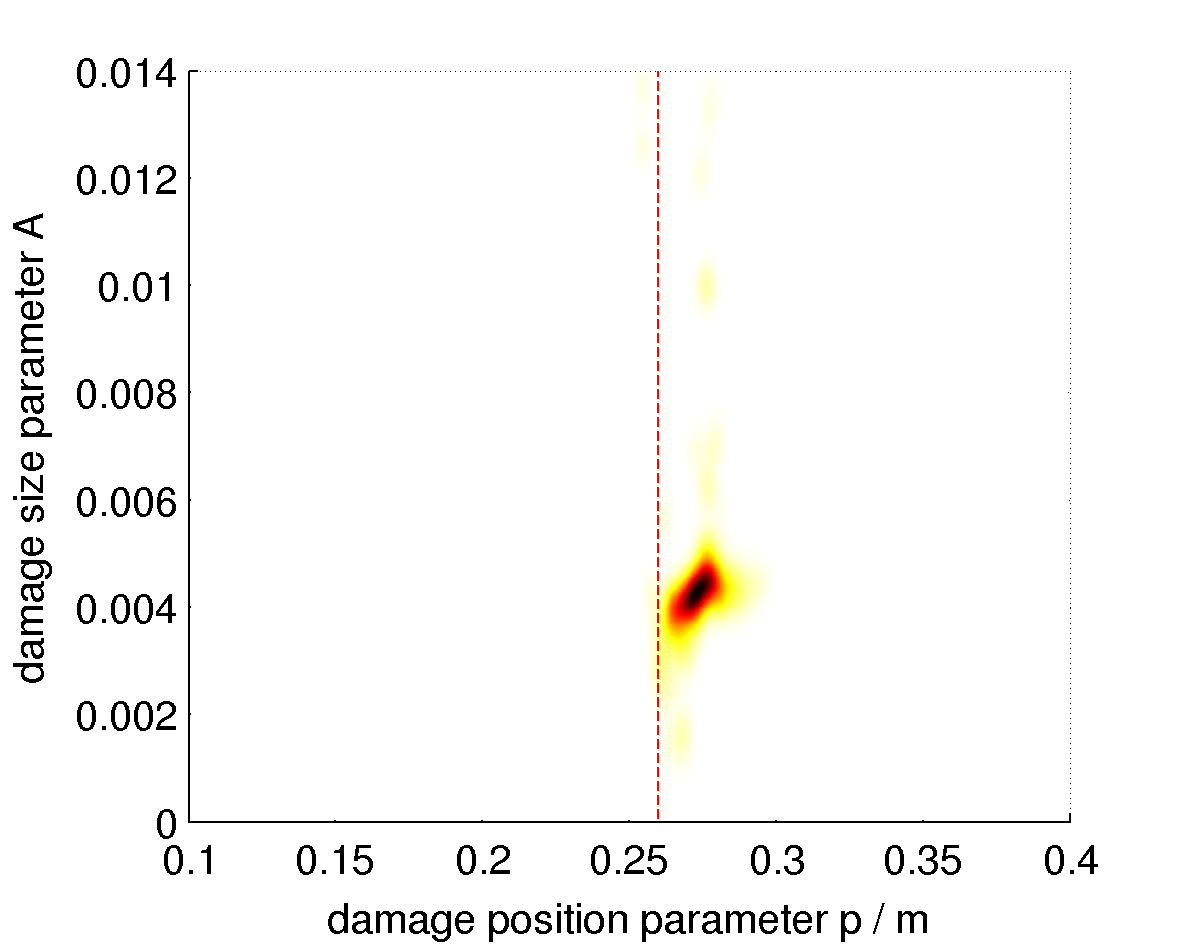}
\includegraphics[scale=0.6]{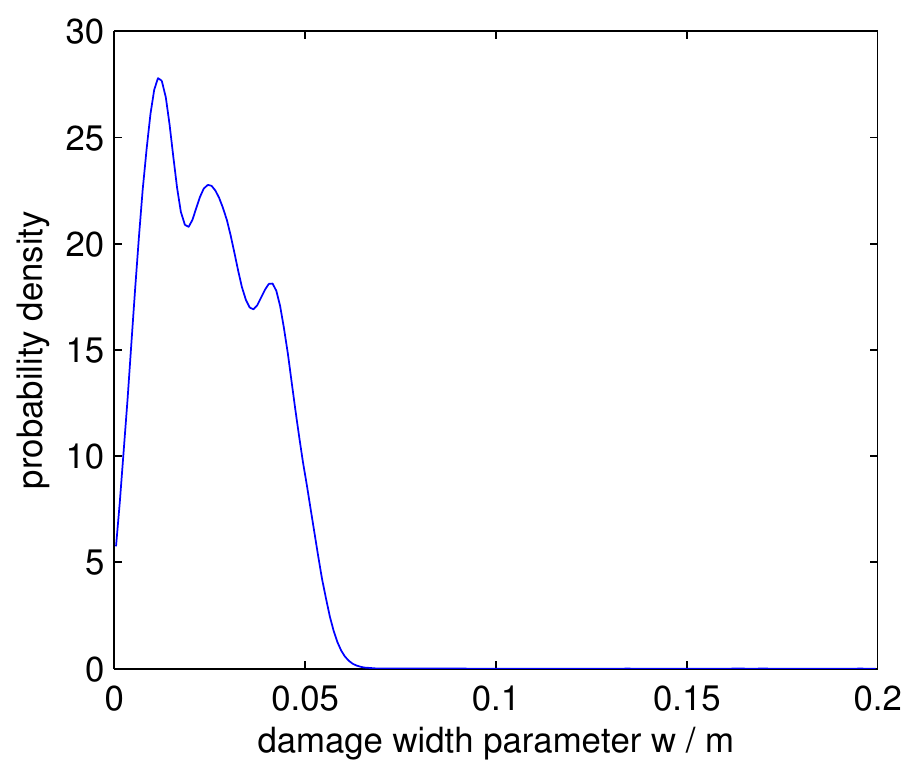}
\end{center}
\caption{Results for experimental damage case 3. Top left: statistics of the EnKF results. Top right: statistics of the regularization method results. Bottom left: marginal density of $(p,A)$ obtained with SACOM. The vertical red line represents the location of the afflicted damage. Bottom right: marginal density of $w$ obtained with SACOM.}\label{f:experData_d_case3}
\end{figure}

In this case the damage was a 15 mm deep and 1 mm wide slot at 260 mm from the fixed end of the beam. Figure \ref{f:experData_d_case3} shows the results from the three methods.

The EnKF now gives a clear indication of damage near the correct location. Damage is also seen near the free end of the beam, which is again interpreted as a known incorrect detection. The regularization method detects damage at the correct location, and an increase in the damage amount is seen. As before, the SACOM results suggest an increase in the amount of damage while the location is near the correct location with high probability. The affected area width remains narrow with high probability.

\subsubsection{Damage case 4}

\begin{figure}[htb]
\begin{center}
\includegraphics[scale=0.6]{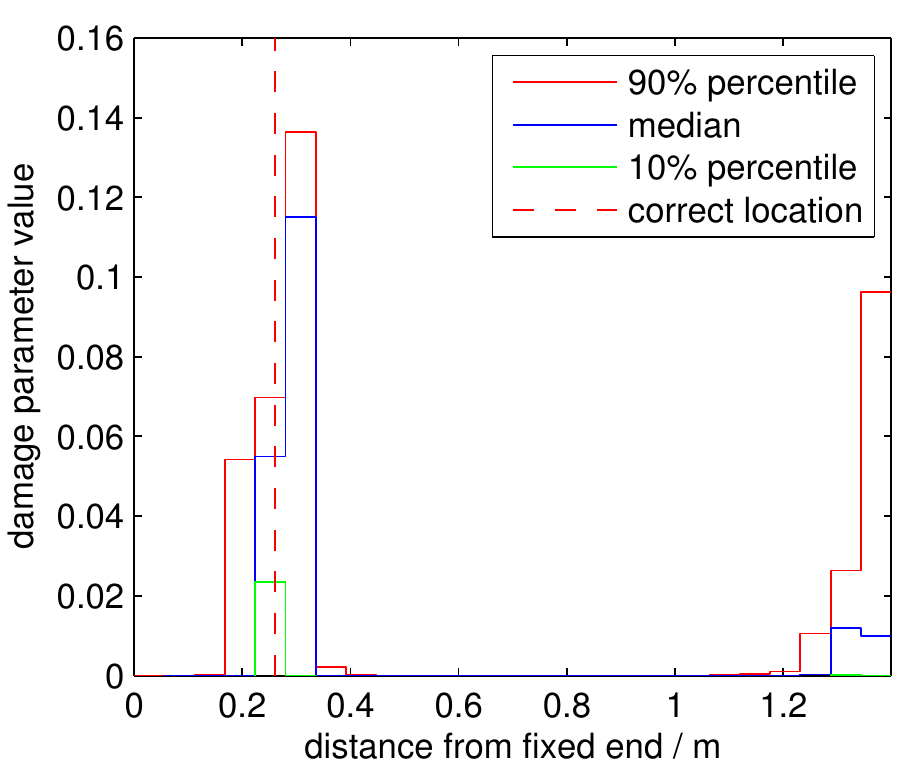}
\includegraphics[scale=0.6]{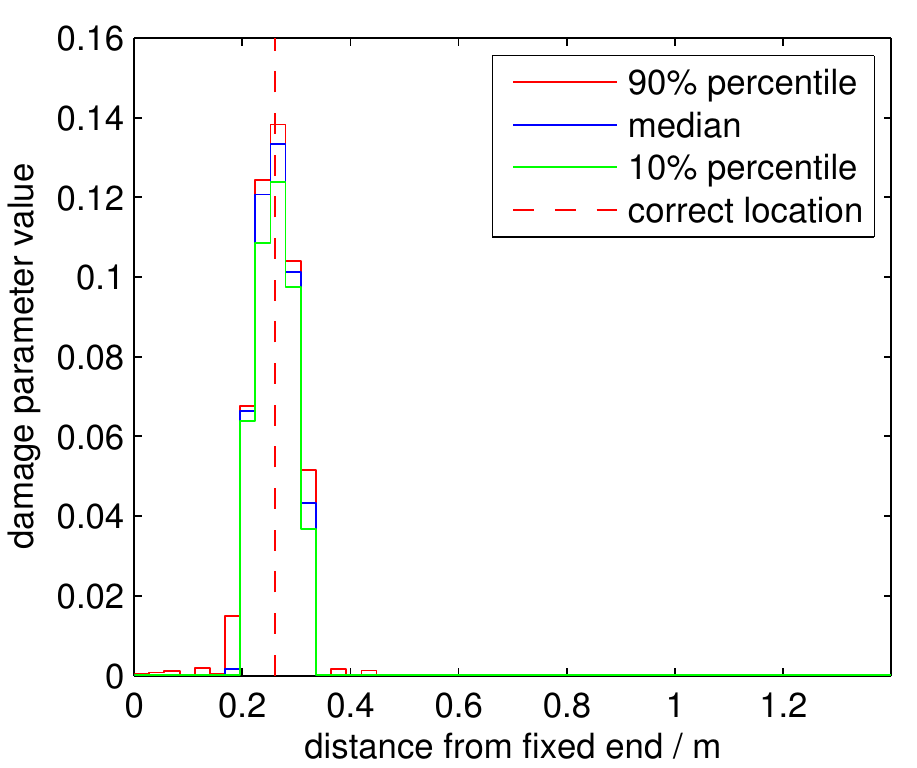} \\
\includegraphics[scale=0.6]{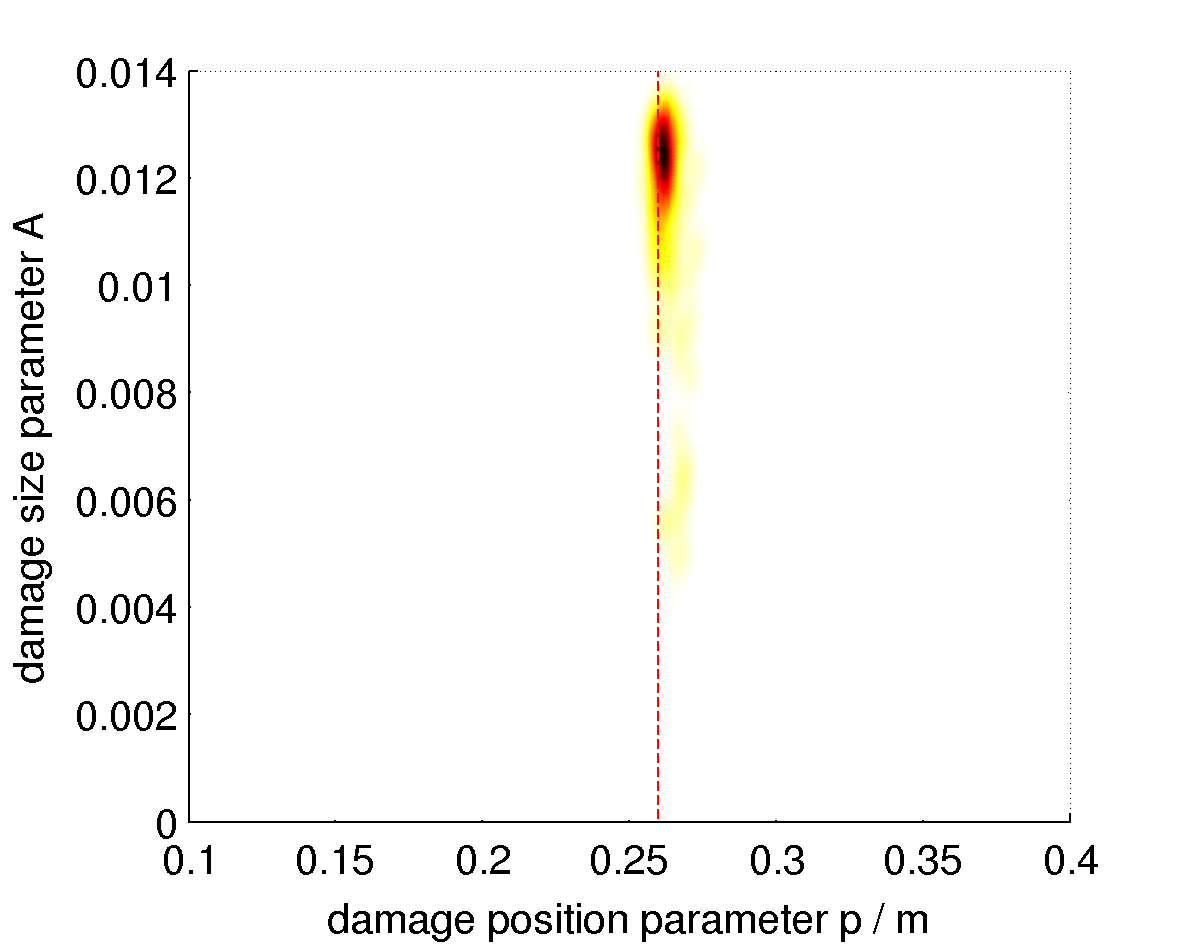}
\includegraphics[scale=0.6]{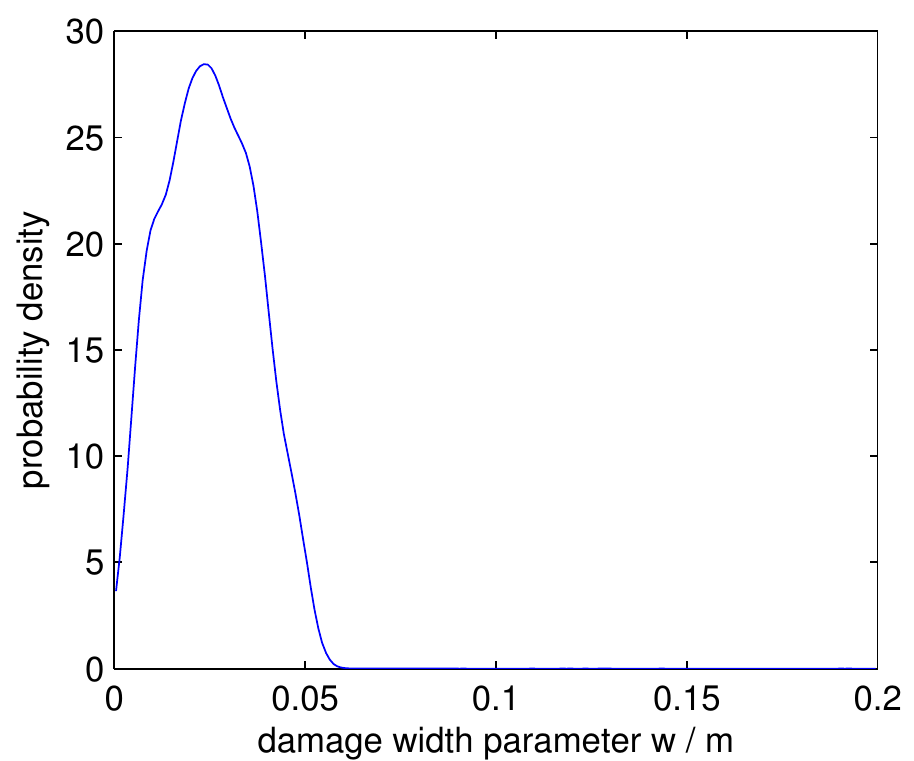}
\end{center}
\caption{Results for experimental damage case 4. Top left: statistics of the EnKF results. Top right: statistics of the regularization method results. Bottom left: marginal density of $(p,A)$ obtained with SACOM. The vertical red line represents the location of the afflicted damage. Bottom right: marginal density of $w$ obtained with SACOM.}\label{f:experData_d_case4}
\end{figure}

In this case the damage was a 20 mm deep and 1 mm wide slot at 260 mm from the fixed end of the beam. Figure~\ref{f:experData_d_case4} shows results of the methods in this damage case.

The EnKF shows a clear indication of damage at the correct location. Damage is still seen near the free end of the beam, which is again interpreted as a known incorrect detection. The regularization method continues to detect damage at the correct location, with the reconstructed damage size having increased from the previous damage case. As with the simulated data for the largest damage case, the SACOM approach produces a density in this case that suggests with high probability the damage is located around the correct location, and while the amplitude of the damage has increased, there is more variability in this value.

% ----------------------------------------------------------------
% ----------------------------------------------------------------
\section{Conclusions}\label{S:Conclusions}
% ----------------------------------------------------------------
% ----------------------------------------------------------------

We have described the frameworks, analysis of data, and demonstrated the utility of interfacing three distinct UQ methodologies for determining and quantifying uncertainty in damage parameters for a vibrating beam. Each method is defined distinctly with respect to the parameter, model, and data spaces along with the physical and numerical maps between these spaces. The EnKF is easiest to implement using ``time-series'' data while the other approaches use modal, i.e.~frequency, data.

The results from the EnKF provide numerical validation for the choice of penalty term in the regularization results. Subsequently, the regularization results suggest a certain parameterized subspace for a damage parameter is appropriate for the SACOM approach. The importance of this can not be understated. The damage parameter may otherwise be modeled as a discretized random field, i.e., a function in an infinite dimensional function space projected onto a discretization. In fact, this is what is done for the EnKF and regularization results. However, in the computational framework we described, applying the SACOM approach on this space, while possible, implies we carry out computational measure theory on a 50-dimensional space. The resulting probability measure is extremely difficult to visualize in such a high dimensional space even though we can identify events of high probability. Moreover, to control the numerical errors in the approximation of such events may require a significant amount of computational resources in terms of a large number of parameter samples and associated model solves if non-adaptive sampling approaches are used \cite{Butler2014b}.

We could further exploit the regularization results by limiting the ranges of the parameters defining the damage parameter for the SACOM method, i.e., we can define small neighborhoods around certain functions suggested by regularization in order to build probability distributions using the SACOM approach. Furthermore, we can vary the penalty terms in the regularization results and use the subsequent results to define more complicated parameterizations of functions and neighborhoods around these functions for computation of probability densities using the SACOM approach. This will be the subject of future work.

%% The Appendices part is started with the command \appendix;
%% appendix sections are then done as normal sections
%% \appendix

%% \section{}
%% \label{}

%% References
%%
%% Following citation commands can be used in the body text:
%% Usage of \cite is as follows:
%%   \cite{key}         ==>>  [#]
%%   \cite[chap. 2]{key} ==>> [#, chap. 2]
%%

%% References with bibTeX database:

%\bibliographystyle{elsarticle-num}
\bibliographystyle{plain}
\bibliography{References}

%% Authors are advised to submit their bibtex database files. They are
%% requested to list a bibtex style file in the manuscript if they do
%% not want to use elsarticle-num.bst.

%% References without bibTeX database:

% \begin{thebibliography}{00}

%% \bibitem must have the following form:
%%   \bibitem{key}...
%%

% \bibitem{}

% \end{thebibliography}

\end{document}